\documentclass[reqno,12pt] {article}
\usepackage{amsmath, amsfonts, amssymb, amsthm,hyperref,floatrow}
\usepackage{graphicx, color,dsfont}
\usepackage{xcolor}
 \hypersetup{
    linktoc=page,
    linkcolor=red,          
    citecolor=blue,        
    filecolor=blue,      
    urlcolor=cyan,
    colorlinks=true           
}

\newtheorem{theorem}{Theorem}
\newtheorem{remark}[theorem]{Remark}
\newtheorem{lemma}[theorem]{Lemma}
\newtheorem{proposition}[theorem]{Proposition}

\newtheorem{corollary}[theorem]{Corollary}
\newtheorem{definition}[theorem]{Definition}

\def \P{ \mathbb P  }
\def \E{ \mathbb E  }

\definecolor{remi}{rgb}{1,0,0}
\usepackage{titlesec}
    \titleformat{\section}[hang]
        {\color{remi}{}\bfseries\filcenter\large}
        {\thesection.}
        {0.4em}
        {}[]

\renewcommand{\tilde}{\widetilde}          
\DeclareMathSymbol{\leqslant}{\mathalpha}{AMSa}{"36} 
\DeclareMathSymbol{\geqslant}{\mathalpha}{AMSa}{"3E} 
\DeclareMathSymbol{\eset}{\mathalpha}{AMSb}{"3F}     
\renewcommand{\leq}{\;\leqslant\;}                   
\renewcommand{\geq}{\;\geqslant\;}                   

\def \one{ {\rm 1}\mkern-4.5mu{\rm I} }

\newcommand{\R}{\mathbb{R}}

\newcommand{\N}{\mathbb{N}}
\newcommand{\Q}{\mathbb{Q}}

\title{Renormalization of Critical Gaussian Multiplicative Chaos and KPZ Relation}

\author{{\sc Bertrand Duplantier\thanks{Institut de Physique Th\'{e}orique, CEA/Saclay F-91191 Gif-sur-Yvette Cedex, France. E-mail: \texttt{bertrand.duplantier@cea.fr}.  Partially supported by  grant  ANR-08-BLAN-0311-CSD5 and by the MISTI MIT-France Seed Fund.}},\, {\sc R\'{e}mi Rhodes\thanks{Universit\'e Paris-Dauphine, Ceremade, UMR 7564. E-mail: \texttt{rhodes@ceremade.dauphine.fr}. Partially supported by grant ANR-11-JCJC.}},\,
{\sc Scott Sheffield\thanks{Department of Mathematics, Massachusetts Institute for Technology, Cambridge, Massachusetts 02139, USA. E-mail: \texttt{sheffield@math.mit.edu}.
Partially supported by NSF grants DMS 064558, OISE 0730136 and DMS 1209044, and by the MISTI MIT-France Seed Fund.}}, \\ and {\sc Vincent Vargas\thanks{Universit\'e Paris-Dauphine, Ceremade, UMR 7564. E-mail: \texttt{vargas@ceremade.dauphine.fr}. Partially supported by grant ANR-11-JCJC.}}
}

\begin{document}
\maketitle

\vspace{-.4in}

 \begin{abstract}
Gaussian Multiplicative Chaos is a way to produce a measure on $\R^d$ (or subdomain of $\R^d$) of the form $e^{\gamma X(x)} dx$, where $X$ is a log-correlated Gaussian field and $\gamma \in [0,\sqrt{2d})$ is a fixed constant.  A renormalization procedure is needed to make this precise, since $X$ oscillates between $-\infty$ and $\infty$ and is not a function in the usual sense.  This procedure yields the zero measure when $\gamma=\sqrt{2d}$.

Two methods have been proposed to produce a non-trivial measure when $\gamma=\sqrt{2d}$.  The first involves taking a derivative at $\gamma=\sqrt{2d}$ (and was studied in an earlier paper by the current authors), while the second involves a modified renormalization scheme.  We show here that the two constructions are equivalent and use this fact to deduce several quantitative properties of the random measure. In particular, we complete the study of the moments of the derivative multiplicative chaos, which allows us to establish the KPZ formula at criticality. The case of two-dimensional (massless or massive) Gaussian free fields is also covered.
 \end{abstract}
%

\normalsize


\section{Introduction}

In the eighties, Kahane \cite{cf:Kah}  developed
a continuous parameter theory of multifractal random measures, called Gaussian multiplicative chaos. His efforts were followed by several authors \cite{allez,bacry,Bar,DRSV,cf:DuSh,Fan,sohier,rhovar,cf:RoVa} coming up with various generalizations at different scales. This family of random fields  has found many applications in various fields of science  ranging from turbulence and mathematical finance to  $2d$-Liouville quantum gravity.

Roughly speaking, a Gaussian multiplicative chaos on $\R^d$ or on a bounded domain of $\R^d$ (with respect to the Lebesgue measure) is a random measure that can formally be written as:
\begin{equation}\label{intro:M}
M(dx)= e^{\gamma X(x)-\frac{\gamma^2}{2}\E[X(x)^2]}\,dx
\end{equation}
where $X$ is a centered Gaussian distribution and $\gamma$ a nonnegative real parameter. The situation of interest is when the field $X$ is log-correlated, that is when
\begin{equation}\label{intro:K}
K(x,y)\stackrel{def}{=}\E[X(x)X(y)]= \ln \frac{1}{|x-y|}+g(x,y)
\end{equation}
 for some bounded continuous function $g$.  In this case, $X$ is a Gaussian random {\em generalized function} (a.k.a.\ {\em distribution}) on $\R^d$ that cannot be defined as an actual function pointwise.  Kahane showed that one can nonetheless give a rigorous definition  to \eqref{intro:M}. Briefly, the idea is to cut-off the singularity  of the kernel \eqref{intro:K} occurring at $x=y$, this procedure being sometimes referred to as ultraviolet cut-off. There are several possible cut-off procedures. The cut-off strategy we will use here is based on a white noise decomposition of the Gaussian distribution $X$. In short, we will assume that have at our disposal a family $(X_t)_t$ of centered smooth enough Gaussian fields with respective covariance kernels $(K_t)_t$ such that for $s,t\geq 0$ 
\begin{equation}\label{WNcov}
\E[X_t(x)X_s(y)]=K_{t\wedge s}(x,y)
\end{equation}
where the kernels $(K_t)_t$ are assumed to be increasing towards $K$ and $K_t(x,x)\simeq t$ as $t\to\infty$. 
This family  $(X_t)_t$  of Gaussian fields  approximates $X$ in the sense that we recover the distribution $X$ when letting $t$ go to $\infty$. Having applied this cut-off, it is now possible to define the approximate measure
\begin{equation}\label{intro:Mn}
M^\gamma_t(dx)= e^{\gamma X_t(x)-\frac{\gamma^2}{2}\E[X_t(x)^2]}\,dx.
\end{equation}
The white noise covariance structure \eqref{WNcov} ensures that the family $(M^\gamma_t)_t$ is a positive martingale, so it converges almost surely. The measure $M^\gamma$ is then understood as the almost sure limit of this martingale. The limiting measure $M^\gamma$ is non trivial if and only if $\gamma^2<2d$  (see \cite{cf:Kah}). Important enough, this construction is rather universal in the sense that the law of the limiting measure $M^\gamma$ is insensitive to the choice of the  cut-off family $(X_t)_t$ made to approximate $X$ (see \cite{cf:Kah,cf:RoVa}). One may for instance smooth up the field $X$ with a mollifying sequence $(\rho_t)_t$  to get a cut-off family $X_t=X\star \rho_t$ by convolution and consider once again the associated measures \eqref{intro:Mn}. Then one will get the same law for the limiting measure.

For $\gamma^2\geq 2d$, the measure $M^{\sqrt{2d}}$ as defined by \eqref{intro:M} thus vanishes, giving rise to the issue of constructing  non trivial objects for $\gamma^2\geq 2d$ in any other possible way. In this paper, we pursue the effort initiated in \cite{DRSV} to understand the critical case  when $\gamma^2=2d$. It is shown in \cite{DRSV} that the natural object at criticality is the {\it derivative Gaussian multiplicative chaos}, which can be formally written as
\begin{equation}\label{intro:M'}
M'(dx)= (\sqrt{2d}\E[X(x)^2]-X(x))e^{\sqrt{2d} X(x)-d\E[X(x)^2]}\,dx.
\end{equation}
It is a positive atomless random measure. It can be rigorously defined via cut-off approximations in the same spirit as in \eqref{intro:Mn}. More precisely, the approximations are obtained by differentiating \eqref{intro:Mn} with respect to $\gamma$ at the value $\gamma^2=2d$, hence the term ``derivative". Nevertheless, it is expected that derivative Gaussian multiplicative chaos at criticality, that is $M'$, can be recovered via a properly renormalized version of \eqref{intro:Mn}. In this paper, inspired by analog results in the case of Galton-Watson processes \cite{heyde,seneta} or branching random walks \cite{AidShi,HuShi}, we prove that
\begin{equation}\label{intro:SH}
\sqrt{t}\,M^{\sqrt{2d}}_t(dx)\to \sqrt{\frac{2}{\pi}}\, M'(dx),\quad \text{ in probability as }t\to \infty.
\end{equation}
 This renormalization procedure is sometimes called Seneta-Heyde renormalization  \cite{AidShi,heyde,seneta}. More precisely, we will first prove \eqref{intro:SH} for covariance kernels \eqref{intro:K} that can be rewritten as 
\begin{equation}\label{introstar}
K(x,y)=\int_1^{+\infty}\frac{k(u(x-y))}{u}\,du
\end{equation}
 for some continuous covariance kernel $k$. Though not covering the whole family of kernels of type \eqref{intro:K}, this family of kernels is quite natural since it possesses nice scaling relations (see \cite{allez}) and are called star-scale invariant kernels.  In this case, the covariance kernel $K_t$ of the approximating field $X_t$ matches
\begin{equation}\label{introstart}
K_t(x,y)=\int_1^{e^t}\frac{k(u(x-y))}{u}\,du.
\end{equation}
Our main standing assumption is the compactness of the kernel $k$, which implies a decorrelation property for the fields $(X_t)_t$, namely that the two process $(X_{t+s}(x)-X_s(x))_{t\geq 0}$ and $(X_{t+s}(w)-X_s(w))_{t\geq 0}$ are independent as soon as $ |x-w|\geq e^{-s}$. So this first step does not allow us to treat the case of long range correlated Gaussian distributions, such as the important class  of two-dimensional free fields.  Anyway, property \eqref{intro:SH} establishes the important fact  
  that the derivative multiplicative chaos   appears as the limit of the natural renormalization of the vanishing martingale \eqref{intro:Mn} for $\gamma^2=2d$. Beyond this unifying perspective, this renormalization approach to criticality is convenient to complete the description of the main properties
of the derivative Gaussian multiplicative chaos initiated in \cite{DRSV}. In particular, the Seneta-Heyde renormalization turns out to be crucial for  applying Kahane's convexity inequalities at criticality. Equipped with this tool, we complete the study of the moments of the derivative multiplicative chaos and compute its power-law spectrum.

 We then treat two-dimensional free fields in a second step. Our main argument is based on  a decomposition of free fields: they can be written  as a sum of two parts, one smooth part that encodes   the long range correlations and another part (not necessarily independent) that  is  logarithmically correlated  at short scales, with a decorrelation property weaker than  in our first case. More precisely, this second part possesses a covariance kernel of the type
\begin{equation}\label{introstarbis}
K(x,y)=\int_1^{+\infty}\frac{k(u,x-y)}{u}\,du
\end{equation} 
for some family $(k(u,\cdot))_u$ of continuous covariance kernels. Hence the structure $k(u,\cdot)$ of correlations at scale $u$ now varies along scales. The corresponding cut-off family $(X_t)_t$  now possesses the following decorrelation property: the two process $(X_{t+s}(x)-X_s(x))_{t\geq 0}$ and $(X_{t+s}(w)-X_s(w))_{t\geq 0}$ are independent as soon as $ |x-w|\geq e^{-s}s^\alpha$ for some $\alpha>0$. This extra polynomial term gives us a sufficient leeway to separate  short from long range correlations. We will prove the Seneta-Heyde renormalization \eqref{intro:SH} for the approximating measure $M^{\sqrt{2d}}_t$ involving a cut-off family for the kernel \eqref{introstarbis} as well, hence  treating the free field case. We stress here that this decomposition argument is not specific to free fields nor to dimension two, and thus may be generalized if need be. 

Studying the case of free fields is motivated by applications in $2d$-Liouville quantum gravity. Recall that the authors in \cite{cf:DuSh} constructed a probabilistic and geometrical framework for two dimensional Liouville quantum gravity and the KPZ equation \cite{cf:KPZ},  based on the two-dimensional Gaussian free field (GFF) (see \cite{cf:Da,DFGZ,DistKa,cf:DuSh,GM,cf:KPZ,Nak} and references therein for physics considerations). It consists in taking $X$ equal to the GFF in \eqref{intro:M} and multiplying the measure by a suitable power of the conformal radius, in the case of Dirichlet boundary conditions in a bounded domain. The resulting measure, called the Liouville measure, is  conformally  invariant   (see \cite[section 2]{cf:DuSh}). In this context, the KPZ relation has been proved rigorously \cite{cf:DuSh,cf:RhoVar} (see also \cite{Benj} in the case of Mandelbrot's multiplicative cascades).  This was done in the standard case of Liouville quantum gravity, namely strictly below the critical value of the GFF coupling constant $\gamma$ in the Liouville conformal factor, i.e, for $\gamma<2$ (the phase transition occurs at $\gamma^2=2d$ in dimension $d$, thus at $\gamma=2$ in two dimensions). Recall that the parameter $\gamma$ is related to the so-called \textit{central charge} $c$ of the conformal field theory coupled to gravity, via the famous KPZ result \cite{cf:KPZ}: $$\gamma=(\sqrt{25-c}-\sqrt{1-c})/\sqrt{6},\,\,\,c\leq 1.$$
The critical case $\gamma=2$ thus corresponds  to the value $c=1$.

Our results extends this picture to the critical case $\gamma=2,c=1$ by  constructing the Liouville measure at criticality on a bounded domain $D\subset \R^2$
\begin{equation}\label{intro:critGFF}
M'(dx)= (2\E[X(x)^2]-X(x))e^{2 X(x)-2\E[X(x)^2]} C(x,D)^2\,dx, 
\end{equation}
where $X$ is a GFF on $D$ (say with Dirichlet boundary conditions) and $C(x,D)$ the conformal radius of $D$ viewed from $x$. This construction is based on a white noise decomposition of the GFF, and is also conformally invariant  in spite of its atypical structure. Also, we complete the description of the moments of  the measure $M'$ and compute its power-law spectrum. These properties lead us to achieve the proof of the KPZ relation at criticality, which strongly relies on the Seneta-Heyde renormalization approach. Let us stress that by construction, this proof is valid in any dimension and for any log-correlated Gaussian field for which one can prove a critical Seneta-Heyde renormalization  of the type \eqref{intro:SH}.

\vspace{2mm}
\noindent {\bf Further remarks about atomic measures at criticality.}

In dimension two, the Liouville quantum gravity measure on a domain $D$ is sometimes interpreted as the image of the intrinsic measure of a random surface $\mathcal M$ under a conformal map that sends $\mathcal M$ to $D$.  This type of ``surface'' is highly singular (not a manifold in the usual sense).  In certain limiting cases where the surface develops singular ``bottlenecks'', one expects the image measure on $D$ to become an atomic measure. In a certain sense (that we will not explain here), constructing these atomic measures requires one to replace $\gamma<2$ by a ``dual value'' $\gamma'>2$ satisfying $\gamma \gamma' = 4$.

It is interesting to consider the analogous atomic measure in the critical case $\gamma = \gamma' = 2$ and to think about what its physical significance might be.  We believe that  both the $\gamma=2$ measures (treated in this paper) and their $\gamma'=2$ ``atomic measure'' variants (see below) have been studied in the physics literature before.  However, when reading this literature about $\gamma=\gamma'=2$ Liouville quantum gravity, it is sometimes difficult to sort out which physical constructions correspond to which mathematical objects.  The remainder of this subsection will give a brief history of these constructions and their relationship to the current work.  This discussion may be safely skipped by the reader without specific background or interest in this area.

The issue of mathematically constructing singular  Liouville measures beyond the  phase  transition, namely for $\gamma>2$, and deriving the corresponding (non-standard dual) KPZ relation has been investigated in \cite{BJRV,Dup:houches,PRL}, giving the first mathematical understanding of the so-called {\it duality} in Liouville quantum gravity  (see \cite{Al,Amb,Das,BDMan,Dur,Jain,Kleb1,Kleb2,Kleb3,Kostov:houches} for an account of physics motivations).  It thus remains to complete the mathematical Liouville quantum gravity picture at criticality, i.e. for $\gamma=2$. From the physics perspective, Liouville quantum gravity at criticality has been investigated in \cite{BKZ,GZ,GrossKleban,GrossM,GubserKleban,KKK,Kleb2,kostov91,kostov92,KS,Parisi,Polch,sugino}. The reader is also referred to \cite{DRSV} for a brief summary about the physics literature on Liouville quantum gravity at criticality.  The Liouville measure at criticality presents an unusual dependence on the Liouville field $\varphi$ (equivalent to $X$ here) of the so-called ``tachyon field'' $T(\varphi)\propto \varphi\, e^{2\varphi}$ \cite{KKK,Kleb2,Polch}.  Its integral over a ``background'' Borelian set $A$ generates the quantum area $\mathcal A =\int_A T(\varphi) dx $, that we can recognize as the formal heuristic expression for the \textit{derivative measure} \eqref{intro:M'}. The possibility at criticality of another tachyon field  of the atypical form  $T(\varphi)\propto e^{2\varphi}$ nevertheless appears in \cite{GubserKleban,KKK,Kleb3}. This form seems to heuristically correspond to a measure of  type \eqref{intro:M} (which actually vanishes for $\gamma=2$). At first sight, our result \eqref{intro:SH} here then seems to suggest that, up to the requested renormalization \eqref{intro:SH} of \eqref{intro:M}, the atypical tachyon field would actually coincide with the usual $\varphi\, e^{2\varphi}$ tachyon field.

However, this atypical tachyon field $e^{2\varphi}$ in Liouville quantum gravity has been associated to another,  non-standard, form of the critical $c=1,\gamma=2$ random surface models. Indeed, the   introduction of higher trace terms   in the action of the $c=1$ matrix model of two-dimensional quantum gravity is known to generate a new critical behavior of the random surface \cite{GubserKleban,KKK,Kleb2,sugino}, with an enhanced \textit{critical} proliferation of spherical bubbles connected one to another by microscopic ``wormholes''.  

 In order to model this non-standard critical theory, it might be necessary to modify the measures introduced here by explicitly introducing ``atoms'' on top of them, using the approach of \cite{BJRV,Dup:houches,PRL} for adding atoms to $\gamma<2$ random measures $M^\gamma$ in the description of   the  dual phase of Liouville quantum gravity.   The ``dual Liouville measure'' corresponding to $\gamma < 2$ involves choosing a Poisson point process from $\eta^{-\alpha-1}d\eta M^\gamma(dx)$, where $\alpha = \gamma^2/4 \in (0,1)$, and letting each point $(\eta,x)$ in this process indicate an atom of size $\eta$ at location $x$.
   When $\gamma =2$ and $\alpha = 1$, we can replace $M^\gamma$ with the derivative measure $M'$ \eqref{intro:M'} (i.e., the limit \eqref{intro:SH}), and use the same construction; in this case (since $\alpha = 1$) the measure a.s.\ assigns infinite mass to each positive-Lebesgue-measure set $A \in \mathcal B(\mathbb R^d)$.  It is nonetheless still well-defined as a measure, and all of its (infinite) mass resides on a countable collection of atoms, each with finite mass.   
    Alternatively, one may use standard L\'evy compensation (intuitively, this amounts to replacing an ``infinite measure'' with an ``infinite measure minus its expectation'', interpreted in such a way that the result is finite) to produce
 a random distribution whose integral against any smooth test function is a.s.\ a finite (signed) value.
  One may expect that this construction yields the continuum random measure associated with the non-standard $c=1,\gamma=2$ Liouville random surface with enhanced bottlenecks, as described in \cite{GubserKleban,KKK,sugino}, thus giving a mathematical interpretation to the (formal)  tachyon field $e^{2\varphi}$ that differs from the renormalized measure \eqref{intro:SH}.

%
\vspace{2mm}
{\bf Organization of the paper.} In Section   \ref{setup}, we introduce the so-called star-scale invariant kernels. We recall the main theorems concerning the construction of the derivative multiplicative chaos in Section \ref{sec:der} as stated in \cite{DRSV}.  In Section \ref{sec:renorm}, we state the Seneta-Heyde renormalization for star-scale invariant kernels and its consequences as moment estimates or power law spectrum of the derivative multiplicative chaos. In Section \ref{sec:KPZ}, we explain how to derive the   KPZ formula at criticality from these moment estimates. Finally we extend these results to Gaussian free fields  in Section \ref{freefieldspres}. The appendix is devoted to the proofs.

\vspace{2mm}
{\bf Acknowledgements: }The authors wish to thank the referees for their careful reading of the manuscript, which helped to improve significantly the final version of this paper.
\section{Setup}\label{setup}
\subsection{Notations}

For a Borelian set $A\subset \R^d$, $\mathcal{B}(A)$ stands for the Borelian sigma-algebra on $A$. All the considered fields are constructed on the same probability space $(\Omega,\mathcal{F},\P)$. We denote by $\E$ the corresponding expectation. Given a Borelian set $A\subset \R^d$, we denote by $A^c$ its complement in $\R^d$. The relation $f\asymp g$ means that there exists a positive constant $c>0$ such that $c^{-1}f(x)\leq g(x)\leq c f(x)$ for all $x$.

\subsection{$\star$-scale invariant kernels}\label{sec:defstar}
Here we introduce the Gaussian fields that we will use in the main part of the paper.  We consider a family of  centered stationary Gaussian processes $((X_{t}(x))_{x \in \R^d})_{t\geq 0}$ where, for each $t\geq 0$, the process $(X_{t}(x))_{x \in \R^d}$ has  covariance given by:
\begin{equation} \label{corrX}
K_t(x)=\E[X_{t}(0)X_{t}(x)]= \int_1^{e^t}\frac{k(ux)}{u}\,du
\end{equation}
for some  covariance kernel $k$ satisfying $k(0)=1$, Lipschitz at $0$ and vanishing outside a compact set. Actually, this condition of compact support is the main input of our assumptions. It may be skipped in some cases: this will be discussed in section \ref{freefieldspres} together with the section of associated proofs \ref{freefields}.

 We also assume that the process $(X_t(x)-X_s(x))_{x\in\R^d}$ is independent of the processes $\big((X_u(x))_{x\in\R^d}\big)_{u\leq s}$ for all $s<t$. In other words, the mapping $t\mapsto X_t(\cdot)$ has independent increments. As a byproduct, for each $x\in\R^d$, the process $t\mapsto X_t(x)$ is a Brownian motion. Such a construction of Gaussian processes is carried out in \cite{allez}.
For $\gamma\geq 0$, we consider the approximate   Gaussian multiplicative chaos $M^\gamma_{t}(dx)$ on $\R^d$:
\begin{equation}\label{defchaos}
M^\gamma_{t}(dx)=e^{\gamma X_{t}(x)-\frac{\gamma^2}{2}\E[X_{t}(x)^2]}dx
\end{equation}

It is well known \cite{allez,cf:Kah} that, almost surely, the family of random measures $(M^\gamma_t)_{t>0}$ weakly converges as $t\to \infty$ towards a random measure $M^\gamma$, which is non-trivial for $\gamma^2<2d$. The purpose of this paper is to investigate the phase transition, that is $\gamma^2=2d$. Recall that we have \cite{DRSV,cf:Kah}:
\begin{proposition}
For $\gamma^2=2d$ (and also for $\gamma^2>2d$), the standard construction \eqref{defchaos} yields a vanishing limiting measure:
\begin{equation} \label{eqn::gamma2limitzero}
 \lim_{t \to \infty} M^{\gamma} _t(dx)=0\quad \text{almost surely}.
 \end{equation}
\end{proposition}
One of the main purposes of this article is to give a non trivial renormalization of the family $(M^{\sqrt{2d}} _t)_t$. We stress that a suitable renormalization should yield a non trivial solution to the lognormal star-equation:
\begin{definition}{\bf Log-normal $\star$-scale invariance.} \label{def_1}
A random Radon measure $M $ is said to be lognormal $\star$-scale invariant if for all $0<\varepsilon\leq 1$, $M $ obeys the cascading rule
\begin{equation}\label{star}
 \big(M (A)\big)_{A\in\mathcal{B}(\R^d)}\stackrel{law}{=} \big(\int_Ae^{ \omega_{\varepsilon}(r)}M^\varepsilon(dr))\big)_{A\in\mathcal{B}(\R^d)}
\end{equation}
where $\omega_{\varepsilon}$ is a stationary stochastically continuous Gaussian process and $M^{ \varepsilon}$ is a random measure independent from $\omega_{\varepsilon}$ satisfying the scaling relation
 \begin{equation}\label{star1}
\big(M^{ \varepsilon}(A)\big)_{A\in\mathcal{B}(\R^d)}\stackrel{law}{=} \big(M (\frac{A}{\varepsilon})\big)_{A\in\mathcal{B}(\R^d)}.
\end{equation} 
The set of relations \eqref{star}+\eqref{star1} will be referred to as $\star$-equation.\qed
\end{definition}

Let us  mention that the authors in \cite{allez} have proved that, for $\gamma^2<2d$, the measure $M^\gamma$ is   lognormal $\star$-scale invariant with
\begin{equation}\label{staromega}
\omega_\varepsilon(r)=\gamma X_{\ln \frac{1}{\varepsilon}}(r)+(\frac{\gamma^2}{2}+d)\ln  \varepsilon,
\end{equation}
where $X_{\ln \frac{1}{\varepsilon}}$ is the Gaussian process introduced in \eqref{corrX}. Furthermore this scaling relation still makes perfect sense when the scaling factor $\omega_\varepsilon $ is given by \eqref{staromega} for the value $\gamma^2=2d$.  Therefore, to define a natural Gaussian multiplicative chaos with respect to a $\star$-scale invariant kernel at the value $\gamma^2=2d$, one has to look for a solution to this equation when the scaling factor is given by \eqref{staromega} with $\gamma^2=2d$ and conversely, each random measure  candidate for being a Gaussian multiplicative chaos with respect to a $\star$-scale invariant kernel at the value $\gamma^2=2d$ must satisfy these relations. In \cite{DRSV}, a non trivial solution has been constructed, namely  the derivative multiplicative chaos. Since it is conjectured that all the non trivial ergodic solutions to this equation (actually we also need to impose a sufficient decay of the covariance kernel of the process $\omega_\varepsilon$, see \cite{allez} for further details) are equal up to a multiplicative factor, it is expected that a non trivial renormalization of the family $(M^{\sqrt{2d}} _t)_t$ converges towards the derivative multiplicative chaos. Proving this is the first purpose of this paper. The second purpose is to prove that the derivative multiplicative chaos satisfies the KPZ formula.

\begin{remark}
As observed in \cite{DRSV}, we stress that the main motivation for considering $\star$-scale invariant kernels is the connection between the associated random measures and the $\star$-equation. Nevertheless, we stress here that  our proofs straightforwardly adapt to the important exact scale invariant kernel $K(x)=\ln_+\frac{1}{|x|}$ constructed in \cite{bacry,rhovar} in dimension $1$ and $2$. Why only $d=1$ or $d=2$? First because this kernel is positive definite only for $d\leq 3$. Second, because a suitable white noise approximation is known only for $d\leq 2$. Variants are constructed in any dimension in \cite{rhovar} but they do not possess any decorrelation property for $d\geq 3$.
\end{remark}

\section{Derivative multiplicative chaos}\label{sec:der}

A way of constructing a solution to the $\star$-equation at the critical value $\gamma^2=2d$ is to introduce the derivative multiplicative chaos $M_{t}'(dx)$ defined by:
\begin{align*}
M_{t}'(dx)& = (\sqrt{2d}\,t-X_{t}(x))e^{\sqrt{2d}X_{t}(x)-d\E[X_{t}(x)^2]}dx .
\end{align*}
It is plain to see that, for each open bounded set $A\subset\R^d$,  the family $(M_{t}'(A))_t$ is a martingale. Nevertheless, it is neither nonnegative nor regular. It is therefore not obvious that such a family converges towards a (non trivial) positive limiting random variable. The following theorem has been proved in \cite{DRSV}:

\begin{theorem}\label{mainderiv}
For each bounded open set $A\subset \R^d$, the martingale $(M'_t(A))_{t\geq 0}$ converges almost surely towards a  positive random variable denoted by $M'(A)$, such that $M'(A)>0$ almost surely. Consequently, almost surely, the (locally signed) random measures $(M'_t(dx))_{t\geq 0}$ converge weakly as $t\to \infty$ towards a positive random measure $M'(dx)$. This limiting measure has full support and is atomless. Furthermore, the measure $M'$ is a solution to the $\star$-equation \eqref{star} with $\gamma=\sqrt{2d}$.
\end{theorem}

\section{Renormalization}\label{sec:renorm}

The main purpose of this paper is to establish that the derivative multiplicative chaos can be seen as the limit of a suitable  renormalization of the family $(M^{\sqrt{2d}}_t)_t$.

\begin{theorem}\label{seneta}
The family  $(\sqrt{ t}M_t^{\sqrt{2d}})_t$ converges in probability as $t\to \infty$ towards a non trivial limit, which turns out to be the same, up to a multiplicative constant, as the limit of the derivative multiplicative chaos. More precisely, we have for all bounded open set $A$:
$$\sqrt{ t}M_t^{\sqrt{2d}}(A)\to \sqrt{\frac{2}{\pi}} M'(A),\quad\text{in probability as }t\to\infty.$$
 \end{theorem}

 The main advantage of this renormalization approach is to make the derivative multiplicative chaos appear as a limit of integrals over exponentials of the field: this is useful to use Kahane's convexity inequality (see \eqref{lem:cvx}). We can then prove:

\begin{corollary}\label{momneg}
The positive random measure $M'(dx)$   possesses moments of order $q$ for all $ q <1$. Furthermore, for all $q>0$ and every non-empty bounded open set $A$, we have
$$\sup_{t \geq 1}\E\Big[\Big(\frac{1}{\sqrt{t}M_t^{\sqrt{2d}}(A)}\Big)^q\Big]<+\infty.$$
\end{corollary}
Actually, only the first part of the above corollary directly results from Kahane's convexity inequalities via comparison to  multiplicative cascades where similar results are proved \cite{HuShi,BK}. The second point concerning the supremum is highly non trivial: though we know that the limit of the family $(\sqrt{t}M_t^{\sqrt{2d}}(A))_t$, namely $M'$, possesses negative moments  of all order, it is not obvious to transfer this statement to the family  $(\sqrt{t}M_t^{\sqrt{2d}}(A))_t$ with a help of a simple conditioning argument, mainly because $(\sqrt{t}M_t^{\sqrt{2d}}(A))_t$ is not a martingale.

We can then determine the  power law spectrum of the  random measure $M'$:
\begin{corollary}\label{prop:pw}
The power law spectrum of the random measure $M'$ is given for $0\leq q <1$ by
\begin{equation}\label{pwbar}
\xi(q)= 2dq-dq^2.
\end{equation}
More precisely, for all $q<1$, we may find a constant $C_q$ such that for any bounded open subset  $A$ of $\R^d$:
$$  \E[M'(\lambda A)^q]\asymp C_q\lambda^{\xi(q)},\quad \text{ as }\lambda\to 0.$$
\end{corollary}

\section{KPZ formula}\label{sec:KPZ}


In this section, we investigate the KPZ formula for the derivative multiplicative chaos. The KPZ formula is a relation between the Hausdorff dimensions of a given set $A$ as measured by the Lebesgue measure or $M'$. So we first recall how to define these dimensions. Given an atomless Radon measure $\mu$ on $\R^d$ and $s\in [0,1]$, we define
$$H^{s,\delta}_\mu(A)= \inf \big\{\sum_k \mu(B_k)^{s} \big\}$$ where the infimum runs over all the coverings $(B_k)_k$ of $A$ with closed Euclidean balls (non necessarily centered at $A$) with radius $r_k\leq \delta$. Clearly, the mapping $\delta>0\mapsto H^{s,\delta}_\mu(A)$ is decreasing. Hence we can define the  s-dimensional $\mu$-Hausdorff measure:
$$H^{s}_\mu(A)=\lim_{\delta\to 0}H^{s,\delta}_\mu(A).$$
The limit exists but may be infinite.  $H^{s}_\mu$ is a metric outer measure on $\R^d$ (see \cite{falc} for definitions). Thus $H^s_\mu$ is a measure on the $\sigma$-field of $H^s_\mu$-measurable sets, which contains all the Borelian sets. The $\mu$-Hausdorff dimension of the set $A$ is then defined as the value
\begin{equation}
{\rm dim}_\mu(A)=\inf\{s\geq 0; \,\,H^s_\mu(A)=0\}.
\end{equation}
Notice that ${\rm dim}_\mu(A)\in [0,1]$.
Since $\mu$ is atomless, the Hausdorff dimension is also characterized by:
\begin{equation}\label{HSD}
{\rm dim}_\mu(A)=\sup\{s\geq 0; \,\,H^s_\mu(A)=+\infty\}.
\end{equation}
This   allows to characterize the Hausdorff dimension as the threshold value at which
the mapping $s\mapsto H^s_\mu(A)$ jumps from $+\infty$ to $0$.

In what follows, given a compact set $K$ of $\R^d$, we define its Hausdorff dimensions
${\rm dim}_{Leb}(K)$ and ${\rm dim}_{M'}(K)$  computed as indicated above with $\mu$ respectively equal to the Lebesgue measure or $M'$, which is possible as they both do not possess atoms. So, a priori, the value of ${\rm dim}_{M'}(K)$ is random. Nevertheless, a straightforward $0-1$ law argument shows that ${\rm dim}_{M'}(K)$ is actually deterministic. We reinforce this intuition by stating:

\begin{theorem}\label{KPZ}{{\bf (KPZ at criticality $\mathbf{\gamma^2=2d}$).}}
Let $K$ be a compact set of $\R^d$. Almost surely, we have the relation
$${\rm dim}_{Leb}(K)=\frac{\xi({\rm dim}_{M'}(K))}{d}=2d\,{\rm dim}_{M'}(K)-d\,{\rm dim}_{M'}(K)^2 .$$
\end{theorem}

\begin{remark}
Let us stress that our proof also allows us to choose $K$ random but independent of the measure $M'$.  \end{remark}

\subsection{Heuristics and open questions on the KPZ formula}
Here, we give a direct  heuristic derivation of Theorem \ref{KPZ} in order to give  a quick intuitive idea of why Theorem \ref{KPZ} is valid, to enlighten the idea behind the proof of Theorem \ref{KPZ}, which involves introducing very particular Frostman measures, and  to develop open questions which can be seen as generalizations or complements to Theorem \ref{KPZ}.

In fact, we will work in the subcritical case $\gamma^2<2d$ (a similar heuristic can be derived for the case $\gamma^2=2d$). Recall that lognormal $\star$-scale invariance for $M^\gamma$, defined by (\ref{intro:M}) with $\gamma^2<2d$, amounts to the following equivalent:
\begin{equation*}
M^\gamma (B(x,r)) \sim r^d e^{\gamma X_{\ln \frac{1}{r}}(x)-\frac{\gamma^2}{2} \ln \frac{1}{r} }
\end{equation*}
where $\sim$ denotes that both quantities are equal up to multiplication by a random factor of order $1$ which does not depend on $r$ (note that the random factor depends on $x$). If we set
$\xi_\gamma(s)=(d+\frac{\gamma^2}{2})s- \frac{\gamma^2}{2} s^2$, it is thus tempting to write the following equivalents for a set $K$:
\begin{align*}
 H^{s}_{M^\gamma}(K) & = \lim_{\delta \to 0}  \inf \big\{\sum_k M^\gamma (B(x_k, r_k))^{s}, \; K \subset \cup_k B(x_k, r_k), \; |r_k| \leq \delta \big\}   \\
& \sim  \lim_{\delta \to 0}  \inf \big\{\sum_k (r_k^d e^{\gamma X_{\ln \frac{1}{r_k}}(x_k)-\frac{\gamma^2}{2} \ln \frac{1}{r_k}})^s, \; K \subset \cup_k B(x_k, r_k), \; |r_k| \leq \delta \big\}    \\
& = \lim_{\delta \to 0}  \inf \big\{\sum_k e^{s \gamma  X_{\ln \frac{1}{r_k}}(x_k)-\frac{s^2 \gamma^2 }{2} \ln \frac{1}{r_k}} r_k^{\xi_\gamma(s)}, \; K \subset \cup_k B(x_k, r_k), \; |r_k| \leq \delta \big\}    \\
& \sim \int_K e^{ s \gamma X(x) - \frac{s^2 \gamma^2 }{2} \E[ X(x)^2]  }   H^{ \xi_\gamma(s)/d }_{Leb}(dx) ,   
\end{align*}
where the last term is a Gaussian multiplicative chaos applied to the Radon measure $H^{ \xi_\gamma(s)/d }_{Leb}(K \cap dx)$ (at least if $H^{ \xi_\gamma(s)/d }_{Leb}(K )< \infty$). These heuristics show that the quantum Hausdorff measure $H^{s}_{M^\gamma}$ should not be too far from $\int_{.} e^{ s \gamma X(x) - \frac{s^2 \gamma^2 }{2} \E[ X(x)^2]  }   H^{ \xi_\gamma(s)/d }_{Leb}(dx)$ (though up to possible logarithmic corrections). In particular, $H^{s}_{M^\gamma}(K)$ is of order $1$ if and only if $ H^{ \xi_\gamma(s)/d }_{Leb}(K)$ is of order $1$. Note that the $e^{ s \gamma X(x) - \frac{s^2 \gamma^2 }{2} \E[ X(x)^2]  }   H^{ \xi_\gamma(s)/d }_{Leb}(dx)$ measure appears in the physics literature on KPZ in the so-called ``gravitational dressing" (see, e.g., \cite{GM}). Naturally, one can ask to what extent the above heuristics can be made rigorous.

\section{Case of two dimensional Free Fields}\label{freefieldspres}
We will show in this section that  analogous results can be proved in the case where $X$ is the GFF on a bounded domain $D \subset \R^2$ or the Massive Free Field (MFF) on   $  \R^2$. Proofs are gathered in section \ref{freefields}.  
\subsection{Gaussian Free Field and Liouville measure at criticality}
Consider a bounded open domain $D$ of $\R^2$. Formally, a GFF on $D$ is a Gaussian distribution with covariance kernel given by the Green function of the Laplacian on $D$ with prescribed boundary conditions (see \cite{She07} for further details). We describe here the case of Dirichlet boundary conditions. The Green function is then given by the formula:
\begin{equation}\label{bounGreen}
G_D (x,y)=  \pi \int_{0}^{\infty}p_D(t,x,y)dt
 \end{equation}
where $p_D$  is the (sub-Markovian) semi-group of a Brownian motion $B$ killed upon touching the boundary of $D$, namely
$$p_D(t,x,y)=P^x(B_{t} \in dy, \; T_D > t)/dy$$ with $T_D=\text{inf} \{t \geq 0, \; B_t\not \in D \}$. The presence of the factor $\pi$ in \eqref{bounGreen}  corresponds to the normalization \eqref{intro:K} for $G_D(x,y)$. The most  direct way to construct a cut-off family of the GFF on $D$ is then to consider  a white noise $W$ distributed on $D\times \R_+$ and to define:
 $$X(x)=\sqrt{\pi}\int_{D\times \R_+}p_D(\frac{s}{2},x,y)W(dy,ds).$$
One can check that $\E[X(x)X(x')]=\pi\int_0^{\infty} p_D(s,x,x')\,ds=G_D(x,x') $. The corresponding cut-off approximations are given by:
$$X_t(x)=\sqrt{\pi}\int_{D\times [e^{-2t},\infty[}p_D(\frac{s}{2},x,y)W(dy,ds).$$
We define the approximating measures
\begin{equation*}
M_t^2(dx)=  e^{2 X_t(x)-2{\rm Var}[X_t(x)^2]}\,dx 
\end{equation*}
and
\begin{equation}\label{expressionmax}
M_t'(dx)= (2{\rm Var}[X_t(x)^2]-X_t(x))e^{2 X_t(x)-2{\rm Var}[X_t(x)^2]}\,dx,
\end{equation}
where ${\rm Var}(X)$ stands for the variance of the random variable $X$. We claim:
\begin{theorem}\label{cvGFF}
The family of random measures $(M_t'(dx))_t$ almost surely weakly converges as $t\to\infty$ towards  a non trivial positive and atomless random measure $M'$, which can also be obtained as
\begin{equation*}
\sqrt{t}\,M^2_t(dx)\to \sqrt{\frac{2}{\pi}}\, M'(dx),\quad \text{ in probability as }t\to \infty.
\end{equation*}
\end{theorem}

Let us now introduce the conformal radius $C(x,D)$ \index{conformal radius} at point $x$ of a planar bounded domain $D$,  defined by
\begin{equation}\label{confrad}
C(x,D):=\frac{1}{|\varphi'(x)|},
\end{equation}
where $\varphi$ is any conformal mapping of $D$ to the unit disc such that 
$\varphi(x)=0$.

\begin{definition}{\bf (Liouville measure at criticality).}
The Liouville measure at criticality is  defined as (we add the superscript $X,D$ to explicitly indicate the dependence on the free field $X$ and the domain $D$)
\begin{equation}\label{def:Lioucrit}
M^{X,D}(dx)=C(x,D)^2M'(dx),
\end{equation}
where  $C(x,D)$ stands for the conformal radius at point $x$ of the planar bounded domain $D$. 
\end{definition}
 
\begin{remark}
The presence of the conformal radius is more easily understood if one wishes to renormalize the family of measures $(e^{2 X_t(x)}\,dx)_t$ with a deterministic renormalizing family. Then we can write
$$ \sqrt{t}e^{-2t}e^{2 X_t(x)}\,dx=\sqrt{t}e^{2{\rm Var}[X_t(x)^2]-2t}\,e^{2X_t(x)-2{\rm Var}[X_t(x)^2]}\,dx.$$
By observing that the family $({\rm Var}[X_t(x)^2]-t)_t$ converges uniformly on the compact subsets of $D$ towards the conformal radius $\ln C(x,D)$ (see \cite[section 2]{cf:DuSh} for circle average approximations and \cite[section 6]{lacoin} in this context) and that the remaining part is nothing but $\sqrt{t}M_t^{2}(dx)$, it is then plain to see that the family $(\sqrt{t}e^{-2t}e^{2 X_t(x)}\,dx)_t$ weakly converges towards $M^{X,D}$, up to the deterministic factor $\sqrt{2/\pi}$.
\end{remark}

One of the most important properties of Liouville quantum gravity is its so-called {\it conformal covariance}, stated as Theorem \ref{conformal} below.  Given another planar domain $\widetilde{D}$ and a conformal map $\psi:\widetilde{D}\to D$, we will denote by $X'$ the (non necessarily centered) GFF on $\widetilde{D}$ defined by
$$\widetilde{X}(x)=X\circ \psi(x)+2\ln|\psi'(x)|.$$
Observe that the family $(X_t(\psi(x)))_t$ is a cut-off approximation of $X\circ \psi $ admitting a white noise decomposition. The measure $M^{\widetilde{X},\widetilde{D}}(dx)$ on $\widetilde{D}$ is then understood as the limit  as $t\to\infty$ of the martingale:
$$M^{\widetilde{X},\widetilde{D}}_t(dx)= \big(2{\rm Var}[X_t(\psi(x))^2]-X_t(\psi(x))-2\ln|\psi'(x)|\big)e^{2X_t(\psi(x))+4\ln|\psi'(x)|-2{\rm Var}[X_t(\psi(x))^2]}  C(x,\widetilde{D})^2 .$$
We claim:
\begin{theorem}{\bf (Conformal invariance). }\label{conformal}
Let $\psi$ be a conformal map from a domain $\widetilde{D}$ to $D$ and consider the GFF $\widetilde{X}$ defined by $\widetilde{X}=X\circ \psi +2\ln |\psi'|$ on  $\widetilde{D}$. Then $M^{X,D}$ is almost surely the image under $\psi$ of the measure $M^{\widetilde{X},\widetilde{D}}$ on $\widetilde{D}$. That is, $M^{\widetilde{X},\widetilde{D}}(A)=M^{X,D}(\psi(A))$ for each Borel measurable $A\subset \widetilde{D}$.
\end{theorem}
The reader is referred to \cite[section 2]{cf:DuSh} for further references and explanations concerning invariance under conformal reparameterization in $2d$-Liouville quantum gravity. 

From Kahane's convexity inequalities and Corollaries \ref{momneg}+\ref{prop:pw}, we deduce:
\begin{corollary}\label{prop:pwGFF}
1) For all set $A\subset D$ such that ${\rm dist}(A,D^c)>0$ and for all $q<1$, the random variable $M^{X,D}(A)$ admits a finite moment of order $q$. \\
2) The power law spectrum of the random measure $M^{X,D}$ is given for $0\leq q <1$ by
\begin{equation}\label{pwbarGFF}
\xi(q)= 4q-2q^2.
\end{equation}
More precisely, for all $q<1$ and all $x\in D$, we may find a constant $C_q$ such that:
$$  \E[M^{X,D}(B(x,r))^q]\asymp C_q r^{\xi(q)},\quad \text{ as }r\to 0.$$
\end{corollary}

We further stress that Kahane's convexity inequalities allows us to show that for all set $A\subset D$ such that ${\rm dist}(A,D^c)>0$ and for all $q<1$, the random variable $M^{X,D}(A)$ admits a moment of order $q$.

Finally we claim:
\begin{theorem}\label{KPZGFF}{{\bf (KPZ at criticality $\mathbf{\gamma =2}$).}}
Let $K$ be a compact set of $D$. Almost surely, we have the relation
$${\rm dim}_{Leb}(K)=4\,{\rm dim}_{M^{X,D}}(K)-2\,{\rm dim}_{M^{X,D}}(K)^2 .$$
\end{theorem}

\subsection{Massive Free Field}
We give here a few explanations in the case of the whole-plane MFF on $\R^2$ with mass $m^2$. The whole-plane MFF is a centered Gaussian distribution with covariance kernel given by the Green function of the operator $2\pi (-\triangle +m^2)^{-1}$ on $\R^2$, i.e. by:
\begin{equation}\label{MFFkernel}
\forall x,y \in \R^2,\quad G_m(x,y)=\int_0^{\infty}e^{-\frac{m^2}{2}u-\frac{|x-y|^2}{2u}}\frac{du}{2 u}.
\end{equation}
The real $m>0$ is called the mass. This kernel can be rewritten as 
\begin{equation} 
G_m(x,y)=\int_{1}^{+\infty}\frac{k_m(u(x-y))}{u}\,du.
\end{equation}
 for some continuous covariance kernel $k_m(z)=\frac{1}{2}\int_0^\infty e^{-\frac{m^2}{2v}|z|^2-\frac{v}{2}}\,dv$.   
$G_m$ is therefore a $\star$-scale invariant kernel as described in subsection
\ref{sec:defstar}, so we stick to the notations of this section. The point is that $k_m$ is not compactly supported.
Yet all the results of this paper apply to the MFF with cut-off approximations given by \eqref{corrX}. Proofs are given in Section \ref{freefields}.

\subsection{Further remarks on the maximum of the discrete GFF}

By analogy with the star-scale invariant case, we conjectured in a first version of this paper that the (properly shifted and normalized) law of the maximum of the discrete GFF (in a domain $D$) on a lattice with mesh $\varepsilon$ going to $0$ converges in law to the sum of $\ln M'(D)$ and an independent Gumbel variable (up to some constant $c$: see also our conjecture 12 in \cite{DRSV}). A few months after this first version, a series of breakthroughs occured on the aforementioned problem. The work \cite{Bram} establishes convergence of the maximum of the discrete GFF in a bounded domain (the framework adopted in this work is that of a discrete GFF with Dirichlet boundary conditions in a domain $D$ on a grid of size $\varepsilon$ going to $0$). Building on this work, the work \cite{Louisdor} obtained the joint convergence in law of the local maxima and the height of these maxima to a Poisson Point Process of the form $Z(dx) \times e^{-2 h}dh$ where $Z$ is a random measure. (Here we adopt the  convention, different from theirs,  that the discrete field at scale $\varepsilon$ has a variance of order $\ln \frac{1}{\varepsilon}$.) Looking at the expression of $Z$ obtained in \cite{Louisdor} (which is written as a limit of the product of $C(x,D)^2$ by an expression of the form \eqref{expressionmax}, where $X_t$ is a regularized GFF obtained by a cut-off procedure different from the one considered in this manuscript), leads to the conclusion that the random measure $Z$ should coincide with the Liouville measure at criticality defined by \eqref{def:Lioucrit} (up to some multiplicative constant). In order to show that that is indeed the case, one must show a universality result similar to \cite{cf:RoVa} in the derivative multiplicative chaos setting, which then  implies universality of the law of the maximum of the GFF on isoradial graphs (see for instance \cite{chelkak} for a definition and useful estimates). The fact that the conformal radius appears in the limit is due to the fact that the field is non stationary and feels the boundary: more precisely, the GFF has smaller fluctuations as one approaches the boundary. One way to get rid of this effect is to consider the following shifted maximum:
\begin{equation*}  
\sup_{x \in A} X_\varepsilon(x)- 2 \ln \frac{1}{\varepsilon} - \ln C(x,D) +\frac{3}{4} \ln \ln \frac{1}{\varepsilon}
\end{equation*}
 where $A$ is some subset of the domain $D$ and $X_\varepsilon$ a discrete GFF at scale $\varepsilon$. In this case, one should recover exactly $M'(dx)$.  

Let us finally mention  another important work \cite{madaule},  where  the case of $\star$-scale invariant kernels is entirely treated, by both establishing convergence and identifying  the limit.


\appendix
\section{Auxiliary results}
We first state the classical ``Kahane's convexity inequalities" (originally written in \cite{cf:Kah} , see also \cite{allez,cf:RoVa} for a proof in English):

\begin{lemma}\label{lem:cvx}
Let $F:\R_+\to \R$ be some convex  function such that, for some positive constants $M,\beta$, $|F(x)|\leq M(1+|x|^\beta)$ for all $  x\in\R_+$. Let $\sigma$ be a Radon measure on the Borelian subsets of $\R^d$. Given a bounded Borelian set $A$, let $(X_r)_{r\in A},(Y_r)_{r\in A}$ be two continuous centered Gaussian processes with continuous covariance kernels $k_X$ and $k_Y$ such that $k_X(u,v)\leq k_Y(u,v)$ for all 
$u,v\in A$.
Then
$$\E\Big[F\Big(\int_A e^{X_r-\frac{1}{2}\E[X_r^2] }\,\sigma(dr)\Big)\Big]\leq \E\Big[F\Big(\int_Ae^{Y_r-\frac{1}{2}\E[Y_r^2] }\,\sigma(dr)\Big)\Big].$$
If we further assume $k_X(u,u)= k_Y(u,u)$ for all $u\in A$,  then for any increasing function $G:\R_+\to \R$:
$$\E\big[G\big(\sup_{x\in A}Y_x\big)\big]\leq \E\big[G\big(\sup_{x\in A}X_x\big)\big].$$
\end{lemma}

\section{Proofs of Section \ref{sec:renorm}}

We denote by $\mathcal{F}_t$ the sigma algebra generated by $\{X_s(x);s\leq t,x\in\R^d\}$ and by $\mathcal{F}$ the sigma algebra generated by $\bigcup_t\mathcal{F}_t$. Given a fixed open bounded set $A\subset\R^d$  and parameters $t,\beta>0$, we  introduce the random variables
\begin{align*}
Z^{\beta}_t(A)=&\int_A(\sqrt{2d}t-X_{t}(x)+\beta )\one_{\{\tau^\beta_x>t\}} e^{\sqrt{2d}X_{t}(x)-dt}\,dx,\quad R^{\beta}_t(A)=&\int_A \one_{\{\tau^\beta_x>t\}} e^{\sqrt{2d}X_{t}(x)-dt}\,dx
 \end{align*}
where, for each $x\in A$, $\tau^\beta_x$ is the $(\mathcal{F}_t)_t$-stopping time   defined by
$$ \tau^\beta_x=\inf\{u>0,X_u(x)-\sqrt{2d}\,u>\beta \}.$$ For $x\in\R^d$, we also define
$$f_t^\beta(x)= (\sqrt{2d}\,t-X_{t}(x)+\beta )\one_{\{\tau^\beta_x>t\}} e^{\sqrt{2d}X_{t}(x)-dt}.$$

It is plain to check that for each $\beta>0$ and each bounded open set $A$, $(Z^{\beta}_t(A))_t$ is a nonnegative martingale with respect to $(\mathcal{F}_t)_t$ such that $\E[Z^{\beta}_t(A)]=\beta |A|$. It is proved in \cite{DRSV} that it is uniformly integrable and therefore almost surely converges towards a non trivial limit.

Recall that, for each $x$ fixed, the process $t\mapsto X_t(x)$ is standard Brownian motion. We will repeatedly use this fact throughout the proof without mentioning it again.

\subsection{Rooted (or Peyri\`ere) measure}\label{rooted}
Since for each $x$, $(f_t^\beta(x))_t$ is a martingale, we can define the probability measure $\Theta^\beta_t$ on $\mathcal{B}(A)\otimes\mathcal{F}_t $ by
\begin{equation}
\Theta^\beta_t=\frac{f_t^\beta(x)}{\beta |A|}dx\,d\P.
\end{equation}
We denote by $ \E_{\Theta^\beta_t}$ the corresponding expectation. In fact, since the above definition defines a pre-measure on the ring $\bigcup_t \mathcal{F}_t$, one can define the rooted measure $\Theta^\beta$ on $\mathcal{B}(A)\otimes\mathcal{F}$ by using Caratheodory's extension theorem. We recover $\Theta^\beta_{| \mathcal{B}(A)\otimes\mathcal{F}_t }=\Theta^\beta_t $. We observe that $\Theta^\beta_t(Z^{\beta}_t(A)>0)=1$ for any $t$.

Similarly, we construct   the probability measure $Q^\beta$ on $\mathcal{F}$ by setting:
$$Q^\beta_{ | \mathcal{F}_t}=\frac{Z_t^\beta(A)}{\beta|A|}\,d\P,$$ which is nothing but the marginal law of $(\omega,x)\mapsto \omega$ with respect to  $\Theta^\beta_t$. Since $(Z_t^\beta(A))_{t \geq 0}$ is a uniformly integrable martingale which converges to a limit $Z^\beta(A)$, we can also define  $Q^\beta$ directly on $\mathcal{F}$ by:
$$Q^\beta=\frac{Z^\beta(A)}{\beta|A|}\,d\P.$$
We state a few elementary results below. In the sequel, when the context is clear, we sometimes identify $\lbrace \emptyset, A  \rbrace \otimes \mathcal{F}_t$ with $\mathcal{F}_t$.
The conditional law of $x$ given $\mathcal{F}_t$ is given by:
$$\Theta^\beta_t(dx |\mathcal{F}_t)=\frac{f^\beta_t(x)}{Z_t^\beta(A)}\,dx.$$
If $Y$ is a $\mathcal{B}(A)\otimes\mathcal{F}_t $-measurable random variable then it has the following conditional expectation given $\mathcal{F}_t$:
\begin{equation}\label{referee1}
\E_{\Theta^\beta_t}[Y | \mathcal{F}_t]=\int_A Y(x,\omega)\frac{f^\beta_t(x)}{Z^\beta_t(A)}\,dx.
\end{equation}
In particular, for any event $E \in \mathcal{F}_t$, we have
$$\E_{Q^\beta}[\one_E\E_{\Theta^\beta_t}[Y|\mathcal{F}_t]]=\E_{Q^\beta}[ \one_E   \int_A Y(x,\omega) \frac{f^\beta_t(x)}{Z^\beta_t(A)}\,dx]= \E_{\Theta^\beta_t}[\one_E Y ].$$

%

 Under $\Theta^\beta_t$, the law of the random process $(\beta+\sqrt{2d} s-X_s)_{s \leq t} $ is that of   a $ 3$-dimensional Bessel process started at $\beta$. This follows first from the Girsanov transform to get rid of the drift term and then from the fact that one is then conditioning the process to stay positive. In what follows, we will use the notation:
\begin{equation}\label{nota:bert}
Y^\beta_s(x)=\beta+\sqrt{2d} s-X_s(x).
\end{equation}
Of course, $Y^0_s(x)$ simply stands for $\sqrt{2d} s-X_s(x)$.

 The proof is inspired from \cite{AidShi}. Nevertheless, we make two remarks. First, most of the auxiliary estimates obtained in  \cite{AidShi} about the minimum of the underlying random walk  are much easier to obtain in our context because of the Gaussian nature of our framework (in particular, the random walk conditioned to stay positive is here a Bessel process). Second, the continuous structure makes correlations  much more intricate to get rid of: we have no spinal decomposition at our disposal, no underlying tree structure, etc. We adapt some arguments developed in \cite{DRSV} at this level.

\begin{remark}
 It is important here to make the distinction between random variable and random function. For instance, for each fixed $z\in\R^d$, the random variable $(x,\omega)\in A\times \Omega \mapsto X_t(z)$ is $\mathcal{F}_t$-measurable. On the other hand, the random variable $(x,\omega) \mapsto X_t(x)$ is not $\mathcal{F}_t$-measurable because of the $x$-dependence. 
\end{remark}

\subsection{Proofs under the Peyri\`ere measure}

The first step is to prove the convergence under $Q^\beta$. This subsection is thus entirely devoted to the proof of the following result:
\begin{proposition}\label{prop:ebeta}
Given $\beta>0$, we have
\begin{align}
\label{renorm:esp}
 \E_{Q^\beta}\Big[\frac{R^{\beta}_t(A)}{Z^{\beta}_t(A)}\Big]= \sqrt{\frac{2}{\pi t}}(1+\varepsilon_1(t)) \\
\label{renorm:var}
 \E_{Q^\beta}\Big[\Big(\frac{R^{\beta}_t(A)}{Z^{\beta}_t(A)}\Big)^2\Big]= \frac{ 2}{\pi t}(1+\varepsilon_2(t)) 
 \end{align}
 for some functions $\varepsilon_1,\varepsilon_2$ going to $ 0$ as $t\to 0$.
As a consequence, under $Q^\beta$, we have
\begin{equation}\label{renorm:cvtheta}
\lim_{t\to \infty}\sqrt{t}\frac{R^{\beta}_t(A)}{Z^{\beta}_t(A)}=\sqrt{\frac{2}{\pi }}\quad \text{ in probability}. \end{equation}
\end{proposition}

\noindent {\it Proof.} It is clear that \eqref{renorm:esp}+\eqref{renorm:var} imply that the variance of the ratio $\frac{\sqrt{t}R^{\beta}_t(A)}{Z^{\beta}_t(A)}$ under $Q^\beta$ goes to $0$ a $t$ goes to $\infty$. Therefore, under $Q^\beta$, this ratio converges in quadratic mean towards $\sqrt{\frac{2}{\pi }}$, and hence in probability, thus proving \eqref{renorm:cvtheta}.

It is plain to establish the relation \eqref{renorm:esp}. First observe that:
\begin{align*}
 \E_{Q^\beta}\Big[\frac{R^{\beta}_t(A)}{Z^{\beta}_t(A)}\Big]&=\frac{1}{\beta |A|}\E[R^{\beta}_t(A)]=\frac{1}{\beta }\E[\one_{\{\tau^\beta_x>t\}}e^{\sqrt{2d}X_t(x)-dt}].
\end{align*}
Let us denote by $B$ a standard one-dimensional Brownian motion. By using the Girsanov transform, we have:
\begin{align*}
 \E_{Q^\beta}\Big[\frac{R^{\beta}_t(A)}{Z^{\beta}_t(A)}\Big]&=\frac{1}{\beta  }\P(\sup_{0\leq u\leq t}B_u\leq \beta) =\frac{1}{\beta  }\P(|B_t|\leq \beta) \simeq \sqrt{\frac{2}{\pi t}}\quad \text{ as }t\to \infty.
\end{align*}
Relation \eqref{renorm:esp} is established.

The proof of \eqref{renorm:var} is much more involved and will be carried out in several steps. To begin with, it may be worth sketching the strategy of the proof:
\begin{enumerate}
\item Observe that (see \eqref{referee1})
\begin{equation}\label{sketch1}
\frac{R_t^\beta(A)}{Z_t^\beta(A)}=\E_{\Theta^\beta_t}\Big[ \frac{1}{Y^\beta_t(x)} |\mathcal{F}_t\Big].
\end{equation}
\item Deduce that
\begin{equation}\label{sketch2}
\E_{Q^\beta}\Big[\Big(\frac{R_t^\beta(A)}{Z_t^\beta(A)}\Big)^2\Big]=  \E_{\Theta^\beta_t}\Big[ \frac{R_t^\beta(A)}{Z_t^\beta(A)}\times \frac{1}{Y^\beta_t(x)}\Big].
\end{equation}
\item The third step consists of subtracting a ball centered at $x$ with radius $e^{-t}$ from the set $A$, call it $B_{x,t}$. This is convenient because if the radius is well chosen the ratio $\frac{R_t^\beta(A\setminus B_{x,t})}{Z_t^\beta(A\setminus B_{x,t})}$ and the weight $\frac{1}{Y^\beta_t(x)}$ will be ``almost" independent. Roughly speaking, the reason why we can subtract a ball is that the measure does not possess atoms. So, at least if the radius of the ball is small enough, it is always possible to subtract a ball centered at $x$ without radically affecting the behaviour of the quantity \eqref{sketch2}. In the forthcoming rigorous proof, we won't base our argument on the fact that the measure is atomless. Instead, we will use estimates on the process $Y^\beta_t(x)$ under the rooted measure, which amounts to the same (as proved in \cite{DRSV}).
\item The last step consists of factorizing \eqref{sketch2} by making use of (almost) independence:
\begin{align*}
\E_{Q^\beta}\Big[\Big(\frac{R_t^\beta(A)}{Z_t^\beta(A)}\Big)^2\Big]&\simeq \E_{\Theta^\beta_t}\Big[\frac{R_t^\beta(A\setminus B_{x,t})}{Z_t^\beta(A\setminus B_{x,t})}\times \frac{1}{Y^\beta_t(x)} \Big]\\
&\simeq \E_{\Theta^\beta_t}\Big[\frac{R_t^\beta(A\setminus B_{x,t})}{Z_t^\beta(A\setminus B_{x,t})} \Big] \E_{\Theta^\beta_t}\Big[\frac{1}{Y^\beta_t(x)} \Big]\\
&\simeq  \E_{Q^\beta}\Big[ \frac{R_t^\beta(A)}{Z_t^\beta(A)}\Big] \E_{Q^\beta}\Big[ \frac{R_t^\beta(A)}{Z_t^\beta(A)}\Big].
\end{align*}
This last quantity is equivalent to $\frac{2}{\pi t}$  for  large $t$. Actually, most of the forthcoming computations are made to justify  that the factorization can be made rigorously. This point is rather  technical and the related computations may appear tedious to the reader.
\end{enumerate}

Now we begin with the rigorous proof of \eqref{renorm:var}. Recall that we use the shorthand \eqref{nota:bert}. Let us first claim:
\begin{lemma} \label{grandO}
 We have
 $$ \E_{Q^\beta}\Big[\Big(\frac{R^{\beta}_t(A)}{Z^{\beta}_t(A)}\Big)^2\Big]=O \Big(\frac{ 1}{  t}\Big)\quad \text{ as }t\to \infty.$$
\end{lemma}

\noindent {\it Proof.}  By Jensen's inequality, we have:
 \begin{align*}
 \E_{Q^\beta}\Big[\Big(\frac{R^{\beta}_t(A)}{Z^{\beta}_t(A)}\Big)^2\Big]&= \E_{Q^\beta}\Big[\Big(\E_{\Theta^\beta_t}\big[\frac{1}{Y^\beta_t(x)}|\mathcal{F}_t\big]\Big)^2\Big] \leq  \E_{Q^\beta}\Big[\E_{\Theta^\beta_t}\big[\frac{1}{Y^\beta_t(x)^2}|\mathcal{F}_t\big]\Big] 
  =  \E_{\Theta^\beta_t}\Big[\frac{1}{Y^\beta_t(x)^2}  \Big]
 \end{align*}
 Since, under   $\Theta^\beta_t$, the law of the process $  (Y^\beta_s(x))_{s \leq t} $ is that of a 3-dimensional Bessel process  starting at $\beta >0$, the lemma follows.\qed

\vspace{1mm}
Now we will decompose the probability space in two parts: a part, call it $E_t$, where we have strong estimates on the process   $  (Y^\beta_s(x))_{s \leq t} $ and an other part that we want to be ``small". More precisely, let $E_t$ be a  $\mathcal{B}(A)\otimes\mathcal{F}_t $-measurable event such that $\E_{\Theta^\beta_t}(\one_{E_t})\to 1$ as $t$ goes to $\infty$. Let
$$\xi_t=\E_{\Theta^\beta_t}\Big[\one_{E_t^c}  \frac{1}{Y^\beta_t(x)} |\mathcal{F}_t\Big].$$
Notice that
\begin{equation}
\frac{R_t^\beta(A)}{Z_t^\beta(A)}=\E_{\Theta^\beta_t}\Big[ \frac{1}{Y^\beta_t(x)} |\mathcal{F}_t\Big]=\xi_t+\E_{\Theta^\beta_t}\Big[ \one_{E_t} \frac{1}{Y^\beta_t(x)} |\mathcal{F}_t\Big].
\end{equation}
We deduce:
\begin{equation}
\E_{Q^\beta}\Big[ \Big(\frac{R_t^\beta(A)}{Z_t^\beta(A)}\Big)^2\Big]=\E_{Q^\beta}\Big[  \frac{R_t^\beta(A)}{Z_t^\beta(A)} \xi_t\Big]+\E_{\Theta^\beta_t}\Big[ \frac{R_t^\beta(A)}{Z_t^\beta(A)} \one_{E_t} \frac{1}{Y^\beta_t(x)} \Big].
\end{equation}
We will treat separately the two terms in the above right-hand side.
The Cauchy-Schwarz inequality and Lemma \ref{grandO} yield
\begin{equation}\label{raj}
\E_{Q^\beta}\Big[  \frac{R_t^\beta(A)}{Z_t^\beta(A)} \xi_t\Big]\leq \E_{Q^\beta}\Big[ \Big(\frac{R_t^\beta(A)}{Z_t^\beta(A)}\Big)^2\Big]^{1/2}\E_{Q^\beta}\Big[ \xi_t^2\Big]^{1/2}\leq \frac{C}{\sqrt{t}}\E_{Q^\beta}\Big[ \xi_t^2\Big]^{1/2}.
\end{equation}
If we can prove that $ \E_{Q^\beta}\big[\xi_t^2 \big]=o\big(\frac{1}{t}\big)$ then \eqref{raj} will tell us that $\E_{Q^\beta}\Big[  \frac{R_t^\beta(A)}{Z_t^\beta(A)} \xi_t\Big]=o\big(\frac{1}{t}\big)$.
Therefore, Proposition \ref{prop:ebeta}, in particular \eqref{renorm:var}, is a consequence of the  following twolemmas:

\begin{lemma}\label{lem:event}
Let $\beta>0$ and $E_t$ be an event such that $\E_{\Theta^\beta_t}(\one_{E_t})\to 1$ as $t$ goes to $\infty$. We have
$$ \E_{Q^\beta}\Big[\xi_t^2 \Big]=o\Big(\frac{1}{t}\Big).$$
\end{lemma}

\begin{lemma}\label{lem:main}
Let $\beta>0$. There exists a family of events $E_t$  such that $\E_{\Theta^\beta_t}(\one_{E_t})\to 1$ as $t$ goes to $\infty$ and
$$\E_{\Theta^\beta_t}\Big[  \frac{R_t^\beta(A)}{Z_t^\beta(A)} \frac{\one_{E_t}}{Y^\beta_t(x)}\Big]\leq  \frac{2}{\pi t} +o\Big(\frac{1}{t}\Big)$$
\end{lemma}

\noindent {\it Proof of Lemma \ref{lem:event}.} Fix $\varepsilon>0$. By Jensen's inequality, we have
\begin{align*}
\E_{Q^\beta}\Big[\xi_t^2 \Big]\leq & \E_{\Theta^\beta_t}\Big[ \one_{E^c_t}\frac{1}{Y^\beta_t(x)^2} \Big]\\
= &  \E_{\Theta^\beta_t}\Big[  \frac{\one_{E_t^c}}{Y^\beta_t(x)^2}\one_{\{\varepsilon\sqrt{t}\leq Y^\beta_t(x)\}}\Big]  +\E_{\Theta^\beta_t}\Big[  \frac{1}{Y^\beta_t(x)^2}\one_{\{\varepsilon\sqrt{t}> Y^\beta_t(x)\}}\Big]\\
= &  \E_{\Theta^\beta_t}\Big[  \frac{\one_{E_t^c}}{Y^\beta_t(x)^2}\one_{\{\varepsilon\sqrt{t}\leq Y^\beta_t(x)\}}\Big] +\frac{1}{\beta}\E \Big[  \frac{\one_{\{\tau^\beta_x>t\}}}{Y^\beta_t(x)}\one_{\{\varepsilon\sqrt{t}> Y^\beta_t(x)\}}e^{\sqrt{2d}X_t(x)-dt }\Big]\\
\leq & \frac{ \E_{\Theta^\beta_t}\big[   \one_{E_t^c}\big]}{\varepsilon^2 t}+\frac{1}{\beta}\E \Big[  \frac{\one_{\{\sup_{[0,t]}X_u(x)\leq \beta\}}}{(\beta-X_t(x)) }\one_{\{X_t(x) > \beta-\varepsilon\sqrt{ t}\}}\Big]
\end{align*}
The third line follows from the translation invariance  of the process $X_t$ and the fourth line from the Girsanov transform. Using the joint law of a Brownian motion together with its maximum, one can check that, for some constant $C$ independent of $t,\varepsilon$, we have:
$$\frac{1}{\beta}\E \Big[  \frac{\one_{\{\sup_{[0,t]}X_u(x)\leq \beta\}}}{(\beta-X_t(x)) }\one_{\{X_t(x) > \beta-\varepsilon\sqrt{ t}\}}\Big]\leq \varepsilon C/t.$$
Since $\E_{\Theta^\beta_t}\big[   \one_{E_t^c}\big]\to 0$ as $t\to\infty$, we deduce:
$$\limsup_{t\to\infty}t\E_{Q^\beta}\Big[\xi_t^2 \Big]\leq \varepsilon C.$$
Since $\varepsilon$ can be chosen arbitrarily small, the proof of the lemma follows.\qed

\vspace{2mm}

\noindent {\it Proof of Lemma \ref{lem:main}.}    We consider a function $h$ such that: 
\begin{equation}\label{defh}
\lim_{t\to \infty}\,\frac{h_t}{t^{\frac{1}{3}}}=+\infty\quad \text{ and }\quad 
\lim_{t\to \infty}\frac{h_t}{t^{\frac{1}{2}}}=0.
\end{equation}
 We define
\begin{align*}
R^{\beta,i}_t(x)&=\int_{B(x,e^{-h_t})\cap A}\one_{\{\tau^\beta_z>t\}}e^{\sqrt{2d}X_t(z)-dt}\,dz &R^{\beta,c}_t(x)&=R^\beta_t(A)-R^{\beta,i}_t(x)  \\
Z^{\beta,i}_t(x)&=\int_{B(x,e^{-h_t})\cap A}f_t^\beta(z)\,dz&Z^{\beta,c}_t(x)&=Z^\beta_t(A)-Z^{\beta,i}_t(x).
\end{align*}
Now we choose the set $E_t$ and define
\begin{align*}
E_t^1&= \{  Z^{\beta,i}_t(x)\leq t^{-2}\} & & E_t^2= \{  R^{\beta,c}_t(x)\leq Z^{\beta,c}_t(x)\}\\
E_t^3&=  \{  R^{\beta,i }_t(x)\leq Z^{\beta,i }_t(x)\} & &
E_t^4= \{h_t^{1/3}\leq Y^\beta_{h_t}(x)\leq h_t\} \\
E_t^5&= \{h_t^{1/3}\leq Y^\beta_{t}(x) \} , & & E_t=E_t^1\cap E_t^2\cap E_t^3\cap E_t^4\cap E_t^5.
\end{align*}

In what follows, $C$ will denote a constant that may change from line to line and that does not depend on relevant quantities. On $E_t$, we have $R^{\beta,i}_t(x)  \leq Z^{\beta,i}_t(x)\leq  t^{-2}$. Therefore
\begin{align*}
\E_{\Theta^\beta_t}\Big[ \frac{R^{\beta,i}_t(x)}{Z_t^\beta(A)}  \frac{1}{Y^\beta_t(x)} \one_{E_t}\Big] & \leq t^{-2} \E_{\Theta^\beta_t}\Big[ \frac{1}{Z_t^\beta(A)}  \frac{1}{Y^\beta_t(x)} \one_{E_t}\Big]  .
\end{align*}
On $E_t$, in particular on $E^5_t$, we have $Y^\beta_t(x)\geq h_t^{1/3}$. Hence
\begin{align*}
\E_{\Theta^\beta_t}\Big[ \frac{R^{\beta,i}_t(x)}{Z_t^\beta(A)}  \frac{1}{Y^\beta_t(x)} \one_{E_t}\Big] & \leq t^{-2} \E_{\Theta^\beta_t}\Big[ \frac{1}{Z_t^\beta(A)}  \frac{1}{h_t^{1/3}}  \Big]  \leq\frac{1}{ t^2 h_t^{1/3} }\E_{Q^\beta}\Big[ \frac{1}{Z_t^\beta(A)} \Big]\leq\frac{1}{\beta h_t^{1/3} |A|t^2} .
\end{align*}
Here we have used $\E_{Q^\beta}\Big[ \frac{1}{Z_t^\beta(A)} \Big]=\frac{1}{\beta|A|}$. Thus we have
\begin{equation}
\E_{\Theta^\beta_t}\big[ \frac{R^{\beta,i}_t(x)}{Z_t^\beta(A)} \frac{1}{Y^\beta_t(x)}\one_{E_t}\big] =o\Big(\frac{1}{t}\Big).
\end{equation}
Let us treat the quantity:
\begin{align*}
\E_{\Theta^\beta_t}&\big[ \frac{R^{\beta,c}_t(x)}{Z_t^\beta(A)} \frac{1}{Y^\beta_t(x)}\one_{E_t}\big]\leq  \E_{\Theta^\beta_t}  \big[ \frac{R^{\beta,c}_t(x)}{Z_t^{\beta,c}(x)} \one_{E^2_t\cap E^4_t}\frac{1}{Y^\beta_t(x)} \big].
\end{align*}
We know that the covariance kernel $k$ appearing in \eqref{corrX} vanishes outside a compact set.  Without loss of generality, we may assume that $k$ vanishes outside the ball $B(0,1)$.  Let us denote by: \\
- $\mathcal{G}_{h_t}$ the sigma-algebra $\mathcal{B}(A)\otimes \mathcal{F}_{h_t}$,\\
- $\mathcal{G}^{out,a}_{h_t}$ the sigma-algebra generated by the random variables $F(x,\omega)$ such that for all $x\in A$, $\omega\mapsto F(x,\omega)$ belongs to the sigma-algebra $\sigma\{X_u(z)-X_{h_t}(z);|z-x|>a, u\geq h_t\}$ (for some $a>0$), \\
- $\mathcal{G}^{in}_{h_t}$ the sigma-algebra generated by the random variables $F(x,\omega)$ such that for all $x\in A$, $\omega\mapsto F(x,\omega)$ belongs to the sigma-algebra $\sigma\{X_u(x)-X_{h_t}(x);u\geq h_t\}$.

Observe that $\mathcal{G}^{out,a}_{h_t}$ and $\mathcal{G}^{in}_{h_t}$ are independent as long as $a>e^{-h_t}$: this comes from the fact that $k$ has compacts support included in the ball $B(0,1)$. We deduce that, conditionally on $\mathcal{G}_{h_t}$,  the random variables $\frac{R^{\beta,c}_t (x)}{Z_t^{\beta,c}(x)} $ and $ \frac{1}{Y^\beta_t(x)} $ are independent since, conditionally on $\mathcal{G}_{h_t}$, $\frac{R^{\beta,c}_t (x)}{Z_t^{\beta,c}(x)}$ is $\mathcal{G}^{out,e^{-h_t}}_{h_t}$-measurable and $Y^\beta_t(x)$ is $\mathcal{G}^{in}_{h_t}$-measurable. Thus we have:
 \begin{align*}
\E_{\Theta^\beta_t} \big[ \frac{R^{\beta,c}_t(x)}{Z_t^\beta(A)} \frac{1}{Y^\beta_t(x)}\one_{E_t}\big]
&\leq \E_{\Theta^\beta_t}\Big[\E_{\Theta^\beta_t}\Big[\frac{R^{\beta,c}_t (x)}{Z_t^{\beta,c}(x)}\one_{ E^2_t\cap E^4_t  } |\mathcal{G}_{h_t}\Big] \E_{\Theta^\beta_t}\Big[  \frac{1}{Y^\beta_t(x)} |\mathcal{G}_{h_t}\Big]\Big].
\end{align*}
To complete the  proof of Lemma \ref{lem:main}, we admit for  a while  the three following lemmas, the proof of which are gathered in the next subsubsection.

\begin{lemma}\label{deadlycomput}
We have
$$ \E_{\Theta^\beta_t}\Big[\frac{1}{Y^\beta_t(x)} |\mathcal{G}_{h_t}\Big] \leq  \sqrt{\frac{2}{\pi (t-h_t)}}   .$$
\end{lemma}

\begin{lemma}\label{deadlycomput2}
We have
$$  \E_{\Theta^\beta_t}\Big[\frac{R^{\beta,c}_t (x)}{Z_t^{\beta, c}(x)} \one_{E^2_t\cap E^4_t} \Big]  \leq \sqrt{\frac{2}{\pi t}}   (1+\varepsilon(t)),$$ for some function $\varepsilon$ such  that   $\lim_{t\to \infty}\varepsilon(t)=0$.
\end{lemma}

\begin{lemma}\label{deadlycomput3}
We have
$$ \lim_{t\to\infty} \Theta^\beta_t(E_t)=1.$$
\end{lemma}

\vspace{2mm}
We conclude the proof of Lemma \ref{lem:main}.  With the help Lemma of \ref{deadlycomput} and \ref{deadlycomput2}, we obtain:
\begin{align*}
\E_{\Theta^\beta_t} \big[ \frac{R^{\beta,c}_t(x)}{Z_t^\beta(A)} \frac{1}{Y^\beta_t(x)}\one_{E_t}\big]\leq & \sqrt{\frac{2}{\pi (t-h_t)}}\E_{\Theta^\beta_t}\Big[\E_{\Theta^\beta_t}\Big[\frac{R^{\beta,c}_t (x)}{Z_t^{\beta,c}(x)}\one_{E^2_t\cap E^4_t} |\mathcal{G}_{h_t}\Big]  \Big]\leq  \frac{2}{\pi t}(1+\varepsilon(t))
\end{align*}
for some function $\varepsilon$ such  that   $\lim_{t\to \infty}\varepsilon(t)=0$.\qed

\subsubsection{Proofs of auxiliary lemmas}

\noindent {\it Proof of Lemma \ref{deadlycomput}.}  By using standard tricks of changes of probability measures, we have:
$$\E_{\Theta^\beta_t}\Big[\frac{1}{Y^\beta_t(x)}|\mathcal{G}_{h_t}\Big] =\frac{1}{f_{h_t}^\beta(x)}\E\Big[ \frac{f^\beta_t(x)}{Y^\beta_t(x)}| \mathcal{G}_{h_t}\Big].$$
Let us compute the latter conditional expectation:
\begin{align*}
\E\Big[  \frac{f^\beta_t(x)}{Y^\beta_t(x)}| \mathcal{G}_{h_t}\Big] 
&=\one_{\{\tau^\beta_x>h_t\}}e^{\sqrt{2d}X_{h_t}(x)-dh_t}\E\Big[ \one_{\{\sup_{[0,t]}X_s(x)-ds\leq \beta\}}e^{\sqrt{2d}(X_{ t}(x)-X_{h_t}(x))-d(t-h_t)}| \mathcal{G}_{h_t}\Big]\\
&=\one_{\{\tau^\beta_x>h_t\}}e^{\sqrt{2d}X_{h_t}(x)-dh_t}F(Y^\beta_{h_t}(x))
\end{align*}
where $F(y)=\E\Big[ \one_{\{\sup_{[0,t-h_t]}X_s(x)-ds\leq y\}}e^{\sqrt{2d}(X_{ t-h_t}(x))-d(t-h_t)} \Big]$. By using the Girsanov transform, we get $F(y)=\E\Big[ \one_{\{\sup_{[0,t-h_t]}X_s(x) \leq y\}} \Big]$.
This quantity is plain to compute since the process $s\mapsto X_s(x)$ is a Brownian motion:
\begin{align*}
F(y)=&\P(|X_{t-h_t}(x)|\leq y)= \frac{1}{\sqrt{2\pi}}\int_{-\frac{y}{\sqrt{t-h_t}}}^{\frac{y}{\sqrt{t-h_t}}}e^{-u^2/2}\,du
\leq  \sqrt{\frac{2}{\pi}}\frac{y}{\sqrt{t-h_t}}.
\end{align*}
Therefore 
$$\E\Big[\frac{f^\beta_t(x)}{Y^\beta_t(x)}| \mathcal{G}_{h_t}\Big]\leq f^\beta_{h_t}(x)\sqrt{\frac{2}{\pi t}} .$$
The lemma follows.\qed

\vspace{2mm}
\noindent {\it Proof of Lemma \ref{deadlycomput2}.}  Let us set $E^6_t=\{h_t^{1/7}\leq \inf_{u\in [h_t,t]}Y_{u}^\beta(x)\}$.  We have:
\begin{align*} 
\E_{\Theta^\beta_t}\Big[  \frac{R^{\beta,c}_t (x)}{Z_t^{\beta,c}(x)} \one_{\{ E^2_t\cap E^4_t\cap E^6_t\}} \Big] 
=  &\E_{\Theta^\beta_t}\Big[\E_{\Theta^\beta_t}\Big[ \frac{R^{\beta,c}_t (x)}{Z_t^{\beta,c}(x)}  \one_{\{ E^2_t\cap E^4_t \}}\one_{\{   E^6_t\}}  |\mathcal{G}_{h_t}\Big]\Big] \\\
= &\E_{\Theta^\beta_t}\Big[ \E_{\Theta^\beta_t}\Big[\frac{R^{\beta,c}_t (x)}{Z_t^{\beta,c}(x)}  \one_{\{ E^2_t\cap E^4_t \}} |\mathcal{G}_{h_t}\Big]\E_{\Theta^\beta_t}[\one_{\{   E^6_t\}}  |\mathcal{G}_{h_t}]\Big] \\
\geq & \E_{\Theta^\beta_t}\Big[ \frac{R^{\beta,c}_t (x)}{Z_t^{\beta,c}(x)}  \one_{\{ E^2_t\cap E^4_t \}}\Big]\inf_{u\in [  h_t^\frac{1}{3}, h_t  ]} \Theta^\beta_t (\inf_{s\in [h_t,t]}Y_{s}^\beta(x)\geq h_t^{1/7} |Y_{h_t}^\beta(x)=u).
\end{align*}
We have used above  the independence of $(Y_s^\beta(x))_{s\leq t}$ and $(R^{\beta,c}_t(x),Z^{\beta,c}_t(x))$,  and thus of   $E^6_t$ and $\frac{R^{\beta,c}_t(x)}{Z^{\beta,c}_t(x)}$, conditionally on $\mathcal{G}_{h_t}$. Now we claim (see the proof below)

\begin{lemma}\label{lem:control}
The following convergence holds:
\begin{equation}\label{control}
\lim_{t\to \infty}\inf_{u\in [h_t^\frac{1}{3} , h_t  ]} \Theta^\beta_t (\inf_{s\in [h_t,t]}Y_{s}^\beta(x)\geq h_t^{1/7} |Y_{h_t}^\beta(x)=u)=1,
\end{equation}
\end{lemma}

Thus we get:
$$  \E_{\Theta^\beta_t}\Big[\frac{R^{\beta,c}_t (x)}{Z_t^{\beta,c}(x)}\one_{E^2_t\cap E^4_t} \Big]  \leq   (1+\rho(t))\E_{\Theta^\beta_t}\Big[  \frac{R^{\beta,c}_t (x)}{Z_t^{\beta,c}(x)} \one_{\{ E^2_t\cap E^4_t\cap E^6_t\}} \Big] ,$$ for some function $\rho$ such  that   $\lim_{t\to \infty}\rho(t)=0$. 
It thus only remains to estimate  the expectation in the right-hand side. Observe that
$$E^2_t\cap E^4_t\cap E^6_t\subset \Big( E^2_t\cap E^6_t \cap  (E^1_t)^c\Big)\cup \Big(E^2_t\cap \{Z_t^{\beta}(A)>1/t\}\cap E^1_t\Big)\cup\Big(E^2_t\cap\{Z_t^{\beta}(A)\leq 1/t\}\Big), $$
in such a way that we get:
\begin{align}\label{plus1}
\E_{\Theta^\beta_t}\Big[  \frac{R^{\beta,c}_t (x)}{Z_t^{\beta,c}(x)} \one_{\{ E^2_t\cap E^4_t\cap E^6_t\}}  \Big] 
\leq & \E_{\Theta^\beta_t}\Big[ \frac{R^{\beta,c}_t (x)}{Z_t^{\beta,c}(x)} \one_{E^2_t\cap E^6_t \cap  (E^1_t)^c}  \Big]  +\E_{\Theta^\beta_t}\Big[ \frac{R^{\beta,c}_t (x)}{Z_t^{\beta,c}(x)} \one_{E^2_t\cap \{Z_t^{\beta}(A)>1/t\}\cap E^1_t}  \Big]\nonumber\\
&+\E_{\Theta^\beta_t}\Big[ \frac{R^{\beta,c}_t (x)}{Z_t^{\beta,c}(x)}\one_{E^2_t\cap \{Z_t^{\beta}(A)\leq 1/t\}} \Big ]\nonumber\\
 \stackrel{def}{=} & A^1_t+A^2_t+A^3_t.
\end{align}
By using Markov's inequality, we have:
\begin{align}\label{plus2}
A^3_t\leq\E_{Q^\beta}[\one_{ \{Z_t^{\beta}(A)\leq 1/t\}}  ]\leq \frac{1}{t}\E_{Q^\beta}[ \frac{1}{Z_t^{\beta}(A)} ]\leq \frac{1}{\beta|A|t} .
\end{align}
Let us now focus on $A^2_t$. On the set  $\{Z_t^{\beta,i}(x) \leq t^{-2}\}\cap\{Z_t^{\beta}(A)>1/t\}$, we have  the estimate $  Z_t^{\beta,i}(x)\leq \frac{1}{t} Z_t^{\beta}(A)$.  Therefore, on this set, we have   $Z_t^{\beta,c}(x)\geq (1-t^{-1})Z_t^{\beta}(A)$. We deduce
\begin{align*} 
 A^2_t&\leq  \frac{1}{1-t^{-1}}\E_{\Theta^\beta_t}\Big[ \frac{R^{\beta,c}_t (x)}{Z_t^{\beta}(A)} \one_{\{ R^{\beta,c}_t(x)\leq Z_t^{\beta,c}(x)\}\cap \{Z_t^{\beta}(A)>1/t\}\cap \{Z_t^{\beta,i}(x) \leq t^{-2}\}}  \Big]  \leq  \frac{1}{1-t^{-1}}  \E_{\Theta^\beta_t}\Big[ \frac{R^{\beta}_t (A)}{Z_t^{\beta}(A)}  \Big]  .
\end{align*}
By using \eqref{renorm:esp}, we deduce
\begin{align}\label{plus3}
A^2_t\leq   \sqrt{\frac{2}{\pi t}}(1+\varepsilon(t)) 
\end{align}
for some function $\varepsilon$ such that $\varepsilon(t)\to 0$ as $t\to \infty$.

It remains to treat $A^1_t$. We have
\begin{align*} 
A^1_t\leq\E_{\Theta^\beta_t}\Big[ \one_{\{ Z_t^{\beta,i}(x) > t^{-2}\}\cap E^6_t}  \Big].
\end{align*}
Let us split $Z_t^{\beta,i}(x)$ in two parts 
$$Z_t^{\beta,i,1}(x) =\int_{B(x,e^{-t})\cap A}f_t^\beta(z)\,dz ,\quad \text{and }\quad Z_t^{\beta,i,2}(x) =\int_{C(x,e^{-t},e^{-h_t})\cap A}f_t^\beta(z)\,dz$$
where $C(x,e^{-t},e^{-h_t})$ stands for the annulus $B(x,e^{-h_t})\setminus B(x,e^{-t})$.
 By using Markov's inequality, we get:
\begin{align*} 
\E_{\Theta^\beta_t}\Big[ \one_{\{ Z_t^{\beta,i}(x) > t^{-2}\}\cap E^6_t}  \Big]\leq &\E_{\Theta^\beta_t}\big[ \one_{\{ 2t^2 Z_t^{\beta,i,1}(x) > 1\} \cap B_t}  \big]+\Theta^\beta_t((B_t)^c)+\E_{\Theta^\beta_t}\big[ \one_{\{ 2t^2 Z_t^{\beta,i,1}(x) > 1\} \cap E^6_t}  \big]  \\
\leq& 2t^2\E_{\Theta^\beta_t}[Z_t^{\beta,i,1}(x)|B_t]+\Theta^\beta_t((B_t)^c)  +2t^2\E_{\Theta^\beta_t}[Z_t^{\beta,i,2}(x)\one_{    E^6_t}  ]  
\end{align*}
where 
$$B_t=\{ t^{1/2-p}\leq Y^\beta_t(x)\leq t^{1/2+p}\} $$ for some $p\in]0,\frac{1}{2}[$ chosen so as to make $\Theta^\beta_t((B_t)^c)=o(t^{-1/2})$. 
This is possible because,  under $\Theta^\beta_t$, the process $Y^\beta(x)$ is a $3$-dimensional Bessel process.
To evaluate $\E_{\Theta^\beta_t}[Z_t^{\beta,i,1}(x)|B_t]$, it is enough to compute
\begin{align}
\E_{\Theta^\beta_t}[&Z_t^{\beta,i,1}(x)|Y^\beta_t(x),x,B_t]=   \int_{B(x,e^{-t})} \E_{\Theta^\beta_t}\Big[f_t^\beta(w)|Y^\beta_t(x),B_t,x\Big]\,dw.
\end{align}   
From the orthogonal projection theorem, observe that
$$X_t(w)=X_t(x)\lambda_t^{x-w}+P^{x,w}_t,$$
where $\lambda_t^{x-w}=\frac{K_t(x-w)}{t}$ and $P^{x,w}_t$ is a centered Gaussian random variable independent of $X_t(x)$ with variance
$$\E[(P^{x,w}_t)^2]=t-\frac{K_t(x-w)^2}{t}.$$
Therefore
\begin{align*}
\E_{\Theta^\beta_t}[&f_t^\beta(w)|Y^\beta_t(x),B_t,x]\\=&    \E_{\Theta^\beta_t}\Big[ (\beta-X_t(w)+\sqrt{2d}t)\one_{\{\tau^\beta_x>t\}}e^{\sqrt{2d}X_t(w)-dt}|Y^\beta_t(x),B_t,x\Big] \\
\leq &    \E_{\Theta^\beta_t}\Big[ |\beta-X_t(x)\lambda_t^{x-w}-P^{x,w}_t+\sqrt{2d}t)|e^{\sqrt{2d}X_t(x)\lambda_t^{x-w}+\sqrt{2d}P^{x,w}_t-dt}|Y^\beta_t(x),B_t,x\Big].
\end{align*}
Now we use the Girsanov transform to get:
\begin{align*}
\E_{\Theta^\beta_t}[&f_t^\beta(w)|Y^\beta_t(x),B_t,x]\\
\leq &    \E_{\Theta^\beta_t}\Big[ \big|\beta-X_t(x)\lambda_t^{x-w}-P^{x,w}_t+\sqrt{2d}t(\lambda_t^{x-w})^2)\big|e^{\sqrt{2d}\lambda_t^{x-w}X_t(x)-dt(\lambda_t^{x-w})^2}|Y^\beta_t(x),B_t,x\Big]\\
=& \E_{\Theta^\beta_t}\Big[ \big|\beta(1-\lambda_t^{x-w})+Y^\beta_t(x)\lambda_t^{x-w}-P^{x,w}_t-\sqrt{2d}t\lambda_t^{x-w}(1-\lambda_t^{x-w})\big|\times\\
&  e^{\sqrt{2d}\lambda_t^{x-w}(\beta +\sqrt{2d}t-Y^\beta_t(x))-dt(\lambda_t^{x-w})^2}|Y^\beta_t(x),B_t,x\Big] .
\end{align*}
To evaluate the above quantity, we first make several observations. By considering $w\in B(x,e^{-t})$, we have $|t-K_t(x-w)|\leq c$ for some constant $c$ that is independent of $t,x,w$. This results from the fact that $k$ is Lipschitz at $0$. We deduce that $1-\lambda_t^{x-w}\leq c/t$ and that $P^{x,w}_t $ has a variance bounded independently of $t,x,w$ such that $w\in B(x,e^{-t})$. Finally, by using the independence of $P^{x,w}_t$ and $X_t(x)$, we deduce that for some constant $C$ independent of $t$ and $|x-w|\leq e^{-t}$:
\begin{align*}
\E_{\Theta^\beta_t}[ f_t^\beta(w)|Y^\beta_t(x),B_t,x]\leq &   C   \E_{\Theta^\beta_t}\Big[ \big(1+Y^\beta_t(x)^2\big) e^{-\sqrt{2d}Y^\beta_t(x)+dt}|Y^\beta_t(x),B_t,x\Big] .
\end{align*}
Therefore:
\begin{align*}
\E_{\Theta^\beta_t}[ Z_t^{\beta,i,1}(x)|Y^\beta_t(x),B_t,x]\leq & C \int_{B(x,e^{-t})}    \E_{\Theta^\beta_t}\Big[ \big(1+Y^\beta_t(x)^2\big) e^{-\sqrt{2d}Y^\beta_t(x)+dt}|Y^\beta_t(x),B_t,x\Big]\,dw\\
\leq &C \int_{B(x,e^{-t})}    \big(1+t^{1+2p}\big) e^{-\sqrt{2d}t^{1/2-p}+dt}\,dw\\
=&C     \big(1+t^{1+2p}\big) e^{-\sqrt{2d}t^{1/2-p} } .
\end{align*}
It is then clear that 
\begin{equation*} 
 t^2\E_{\Theta^\beta_t}[Z_t^{\beta,i,1}(x)|B_t]=o(t^{-1/2}).
 \end{equation*}

To evaluate  $\E_{\Theta^\beta_t}[Z_t^{\beta,i,2}(x)\one_{E^6_t}]$, the strategy is similar. First observe that
\begin{align*}
\E_{\Theta^\beta_t}[Z_t^{\beta,i,2}(x)\one_{E^6_t}]=&\frac{1}{\beta |A|}\int_{x\in A}\Big(\int_{C(x,e^{-t},e^{-h_t})\cap A}\E[\one_{E^6_t}f_t^\beta(w)f^\beta_t(x)]\,dw\Big)\,dx\\
\leq & \frac{1}{\beta |A|}\int_{x\in A}\Big(\int_{C(x,e^{-t},e^{-h_t})\cap A}\E[\one_{\{Y^\beta_{\ln\frac{1}{|x-w|}}(x)\geq h_t^{1/7}\}}f_t^\beta(w)f^\beta_t(x)]\,dw\Big)\,dx.
\end{align*}
We then apply the stopping time theorem. For this, recall that $k$ is assumed to have support included in the ball $B(0,1)$ so that, setting $s_0=\ln\frac{1}{|x-w|}$, the processes $((X_t(x)-X_{s_0}(x))_{t\geq s_0}$ and $((X_t(w)-X_{s_0}(w))_{t\geq s_0}$ are independent. Hence:
\begin{align*}
\E_{\Theta^\beta_t}&[Z_t^{\beta,i,2}(x)\one_{E^6_t}]  \\
\leq & \frac{1}{\beta |A|}\int_{x\in A}\Big(\int_{C(x,e^{-t},e^{-h_t})\cap A}\E\Big[\one_{\{Y^\beta_{s_0}(x)\geq h_t^{1/7}\}}f_{s_0}^\beta(w)f^\beta_{s_0}(x)\Big]\,dw\Big)\,dx.
\end{align*}
Then we can reproduce the previous  computations we made for $Z^{\beta,i,1}(x)$ to get, for some constant $C$ independent of $t$:
 \begin{align*}
\E_{\Theta^\beta_t}&[Z_t^{\beta,i,2}(x)\one_{E^6_t}] \nonumber\\ \leq &   \frac{C}{\beta |A|}\int_{ A} \int_{C(x,e^{-t},e^{-h_t})}  \E\Big[ \one_{\{Y^\beta_{s_0}(x)\geq h_t^{1/7}\}} \big(1+Y^\beta_{s_0}(x)^2\big) e^{-\sqrt{2d}Y^\beta_{s_0}(x)+d{s_0}}f^\beta_{s_0}(x) \Big]\,dwdx \\
\leq& \frac{Ce^{-\sqrt{2d}h_t^{1/7}}}{\beta |A|}\int_{ A} \int_{C(x,e^{-t},e^{-h_t})}  \frac{1}{|w-x|^d}\E\Big[  \big(1+Y^\beta_{\ln\frac{1}{|w-x|}}(x)^2\big) f^\beta_{\ln\frac{1}{|w-x|}}(x) \Big]\,dwdx .
\end{align*}
Thanks to the usual $3d$-Bessel argument, the last expectation is easily seen to be less than $D(1+\ln\frac{1}{|w-x|})$ for some constant $D$ which does not depend on $w,x$. Therefore, we have (for some irrelevant constant $C'$, which may change along lines):
 \begin{align*}
\E_{\Theta^\beta_t} [Z_t^{\beta,i,2}(x)\one_{E^6_t}]  
\leq& \frac{C'e^{-\sqrt{2d}h_t^{1/7}}}{  |A|}\int_{ A} \int_{C(x,e^{-t},e^{-h_t})}  \frac{1}{|w-x|^d}(1+\ln\frac{1}{|w-x|})\,dwdx \\
\leq & C'e^{-\sqrt{2d}h_t^{1/7}} \int_{e^{-t}}^{e^{-h_t} }\frac{1}{r}  (1+\ln\frac{1}{r}  )   \,dr \\
\leq &        C'e^{-\sqrt{2d}h_t^{1/7}} t^2.
\end{align*}
Because of   \eqref{defh}, we deduce $t^2\E_{\Theta^\beta_t}[  Z^{\beta,i,2}(x) \one_{E^6_t}]  \leq  C t^4 e^{-\sqrt{2d}t^{1/21}}  $
in such a way  that
\begin{equation*} 
 t^2\E_{\Theta^\beta_t}[Z_t^{\beta,i,2}(x)\one_{E^6_t}]=o(t^{-1/2}).
 \end{equation*}
Finally, we conclude that
\begin{equation} \label{plus4}
A^1_t=o(t^{-1/2}).
 \end{equation} 
By gathering \eqref{plus1}+\eqref{plus2}+\eqref{plus3}+\eqref{plus4}, the proof is complete.\qed

\vspace{2mm}
\noindent {\it Proof of Lemma \ref{lem:control}.} Recall that, under $\P_{\Theta^\beta}$, the process $(Y^\beta_t(x))_t$ has the law of a $3$-dimensional process. Via coupling arguments, it is plain to check that the infimum is achieved at $u=h_t^{1/3}$. By scaling, we have
\begin{align*} 
\Theta^\beta[ \inf_{s\in [h_t,t]} Y_{s}^\beta(x)\geq h_t^{1/7} |Y_{h_t}^\beta(x)=h_t^{1/3}]=&\Theta^\beta[ \inf_{s\in [0,t-h_t]}Y_{s}^\beta(x)\geq h_t^{1/7} |Y_{0}^\beta(x)=h_t^{1/3}]\\
= &\Theta^\beta\Big[ \inf_{s\in [0,\frac{t-h_t}{h_t^{2/3}}]}Y_{s}^\beta(x)\geq h_t^{1/7-1/3} |Y_{0}^\beta(x)=1\Big]\\
\geq & \Theta^\beta\Big[ \inf_{s\in [0,+\infty[}Y_{s}^\beta(x)\geq h_t^{1/7-1/3} |Y_{0}^\beta(x)=1\Big].
\end{align*}
By taking the limit as $t\to\infty$, we deduce
$$\lim_{t\to \infty}\inf_{u\in [h_t^\frac{1}{3} , h_t  ]}\Theta^\beta_t[\inf_{s\in [h_t,t]}Y_{s}^\beta(x)\geq h_t^{1/7} |Y_{h_t}^\beta(x)=u]\geq \Theta^\beta\Big[ \inf_{s\in [0,+\infty[}Y_{s}^\beta(x)>0|Y_{0}^\beta(x)=1\Big]=1.$$
Indeed, this last quantity is the same as the probability that a $3d$-Brownian motion started at $0$ never hits the point  $(1,0,0)\in\R^3$, which is $1$. 
\qed

\vspace{2mm}
\noindent {\it Proof of Lemma \ref{deadlycomput3}.}  Observe first that we have shown in the proof of Lemma \ref{deadlycomput2} that 
\begin{equation} \label{eqnew1}
 \lim_{t\to\infty} \Theta^\beta_t(E^1_t)=1.
 \end{equation} 
 
Under $\Theta^\beta_t$, the process $(Y^\beta_s(x))_{s \leq t}$ is a $3$-dimensional Bessel process. It is plain to deduce that
\begin{equation}\label{eqnew2}
\lim_{t\to\infty} \Theta^\beta_t (E^4_t)=1\quad \text{and }\quad \lim_{t\to\infty} \Theta^\beta_t (E^5_t)=1.
\end{equation}

Let us prove that 
\begin{equation}\label{eqnew3}
\lim_{t\to\infty} \Theta^\beta_t (E^3_t)=1.
\end{equation}
We define the random variables
$$S_t= \sup_{x\in A}X_t(x)-\sqrt{2d}t+\frac{1}{4\sqrt{2d}}\ln (t+1),\quad\quad S=\sup_{t\geq 0}S_t,$$
and the event
$$B=\big\{S<+\infty\big\}.$$
It is proved in \cite{DRSV} that $\P(B)=1$.
It is plain to deduce that $\Theta^\beta(B)=1$. Then we have
\begin{align*}
Z_t^{\beta,c}(x)=&\int_{A\setminus B(x,e^{-h_t})}Y^\beta_t(z)\one_{\{\tau^\beta_z>t\}}e^{\sqrt{2d}X_t(z)-dt}\,dz\\
\geq  &\int_{A\setminus B(x,e^{-h_t})}\big(\beta+\sqrt{2d}t-\sup_{x\in A}X_t(x)\big)\one_{\{\tau^\beta_z>t\}}e^{\sqrt{2d}X_t(z)-dt}\,dz\\
=& \Big(\beta+\frac{1}{4\sqrt{2d}}\ln (t+1)-S_t\Big)R_t^{\beta,c}(x)\\
\geq& \Big(\beta+\frac{1}{4\sqrt{2d}}\ln (t+1)-S\Big)R_t^{\beta,c}(x).
\end{align*}
Therefore
$$\{\beta+\frac{1}{4\sqrt{2d}}\ln (t+1)-S\geq 1\}\subset \{Z^{\beta,c}_t(x)\geq R^{\beta,c}_t(x)\},$$
from which we deduce $\lim_{t\to\infty} \Theta^\beta_t \big(Z^{\beta,c}_t(x)\geq R^{\beta,c}_t(x)\big)=1$. With the same argument, we prove 
\begin{equation}\label{eqnew4}
\lim_{t\to\infty} \Theta^\beta_t (E^2_t)=1.
\end{equation} 
From \eqref{eqnew1}+\eqref{eqnew2}+\eqref{eqnew3}+\eqref{eqnew4}, we obtain $\lim_{t\to\infty} \Theta^\beta_t (E_t)=1$, thus thus completing the proof of the lemma.\qed

\subsection{Proof of Theorem \ref{seneta}}
In the previous subsection, we have proved the convergence of $\sqrt{t}\frac{R^\beta_t(A)}{Z_t^\beta(A)}$  in probability towards $\sqrt{\frac{2}{\pi}}$ under the measure $Q^\beta$ for any $\beta>0$. This is the content of Proposition \ref{prop:ebeta}. Our objective is now to use this convergence under $Q^\beta$ to establish the convergence under the original probability measure $\P$. Proposition \ref{prop:ebeta} ensures that, for any $\varepsilon>0$,
$$\Q^\beta\Big(\Big|\sqrt{t}\frac{R^\beta_t(A)}{Z_t^\beta(A)}-\sqrt{\frac{2}{\pi}}\Big|>\varepsilon\Big)\to 0,\quad \text{as }t\to\infty.$$
Equivalently
$$\E\Big(Z_t^\beta(A)\one_{\big\{\big|\sqrt{t}\frac{R^\beta_t(A)}{Z_t^\beta(A)}-\sqrt{\frac{2}{\pi}}\big|> \varepsilon\big\}}\Big)\to 0,\quad \text{as }t\to\infty.$$
From \cite{DRSV}, we know that $\sup_t\max_{x\in A}X_t(x)-\sqrt{2d}\,t<\infty$  almost surely. By setting
$$E_R=\{\sup_t\max_{x\in A}X_t(x)-\sqrt{2d}\,t<R\},$$
we obtain an increasing family such that $\P\Big(\bigcup_{R>0} E_R\Big)=1$. From the convergence ($Z_t^\beta(A)$ is nonnegative)
$$\E\Big(Z_t^\beta(A)\one_{\big\{\big|\sqrt{t}\frac{R^\beta_t(A)}{Z_t^\beta(A)}-\sqrt{\frac{2}{\pi}}\big|> \varepsilon\big\}}\one_{E_R}\Big)\to 0,\quad \text{as }t\to\infty,$$
we deduce the convergence of
\begin{equation}\label{cvR}
Z_t^\beta(A)\one_{\big\{\big|\sqrt{t}\frac{R^\beta_t(A)}{Z_t^\beta(A)}-\sqrt{\frac{2}{\pi}}\big|> \varepsilon\big\}}\one_{E_R}\to 0
\end{equation} in probability as $t\to \infty$.
Fix $\beta>R$. On $E_R$ we have $R^\beta_t(A)=M_t^{\sqrt{2d}}(A)$. Concerning $Z_t^\beta(A)$, we observe that, for $\beta>R$, we have
$$\forall t>0,\quad   \beta M_t^{\sqrt{2d}}(A)+M'_t(A)=Z_t^\beta(A). $$
Therefore $\lim_{t\to\infty}Z_t^\beta(A)=M'(A)>0$ on $E_R$ for $\beta>R$ (recall that $M_t^{\sqrt{2d}}(A)\to 0$ as $t\to \infty$). From \eqref{cvR}, we deduce that, necessarily on $E_R$:
$$ \one_{\big\{\big|\sqrt{t}\frac{R^\beta_t(A)}{Z_t^\beta(A)}-\sqrt{\frac{2}{\pi}}\big|> \varepsilon\big\}} \to 0,\quad \text{as }t\to \infty$$
that is
$$ \one_{\big\{\big|\sqrt{t}\frac{M_t^{\sqrt{2d}}(A)}{M'_t (A)}-\sqrt{\frac{2}{\pi}}\big|> \varepsilon\big\}} \to 0,\quad \text{as }t\to \infty.$$
The proof of Theorem \ref{seneta} is over.\qed

\subsection{Other proofs of Section   \ref{sec:renorm}}
We use the comparison with multiplicative cascades set out in the appendix of \cite{DRSV}. The idea is to compare moments of discrete lognormal multiplicative cascades to moments of the family $(M_t)_t$ thanks to Lemma \ref{lem:cvx}.
More precisely and sticking to the notations in \cite{DRSV}, we have
\begin{equation}\label{comp:ccv}
 \E\big[\big(e^{Z} \sqrt{n} M^{\sqrt{2d}}_{n\ln 2}([0,1]^d)\big)^\alpha\big]\leq  \E\big[\big( \sqrt{n} \int_{T}e^{\overline{X}_n(t)-\frac{1}{2}\E[\overline{X}_n(t)^2]}\,\sigma(dt)  \big)^\alpha\big],
\end{equation}
where:\\
$\bullet$ $Z$ is a Gaussian random variable with fixed mean and variance (thus independent of $n$), and independent of the family $(M_t)_t$,\\
$\bullet$ the parameter $\alpha$ belongs to $]0,1[$, making the mapping $x\mapsto x^\alpha$ concave,\\
$\bullet$ $\int_{T}e^{\overline{X}_n(t)-\frac{1}{2}\E[\overline{X}_n(t)^2]}\,\sigma(dt) $ stands for a lognormal multiplicative cascade at generation $n$, the parameters of which are adjusted to be in the critical situation (see \cite{DRSV} for a precise definition).

It is now well established that the right-hand side of \eqref{comp:ccv}  is bounded uniformly with respect to $n$. The reader may consult   \cite{AidShi,BK,HuShi} for instance about this topic.  By Fatou's lemma and Theorem \ref{seneta}, we have
$$ \E\big[\big(  M'([0,1]^d)\big)^\alpha\big] \leq \liminf_{n\to\infty}\E\big[\big( \sqrt{n} M^{\sqrt{2d}}_{n\ln 2}([0,1]^d)\big)^\alpha\big]$$
This shows that the measure $M'$ possesses moments of order $q$ for $0\leq q <1$. We already know that it possesses moments of negative order \cite{DRSV}.

It remains to prove that
$$\sup_{t \geq 1}\E\Big[\Big(\frac{1}{\sqrt{t}M_t^{\sqrt{2d}}(A)}\Big)^q\Big]<+\infty.$$ We write the proof in dimension $d=1$ (generalization to higher dimensions is straightforward) and we assume that the kernel $k$ of \eqref{corrX} vanishes outside the interval $[-1,1]$. Our proof is based on an argument in \cite{Mol}, which we adapt here to get bounds that are uniform in $t$. We work with the ball $A=[0,1]$. The first step consists in writing an appropriate decomposition for the measure $\sqrt{t}M_t^{\sqrt{2}}([0,1])$. We have:
\begin{align*}
M_{t+\ln 8}^{\sqrt{2}}([0,1]) =&\int_0^1e^{\sqrt{2}X_{t+\ln 8}(x)-(t+\ln 8)}\,dx\geq  \sum_{k=0}^3\int_{\frac{k}{4}}^{\frac{2k+1}{8}}e^{\sqrt{2}X_{t+\ln 8}(x)-(t+\ln 8)}\,dx\\
\geq &\sum_{k=0}^3\inf_{x\in [\frac{k}{4},\frac{2k+1}{8}]}e^{\sqrt{2}X_{\ln 8}(x)-\ln 8}\int_\frac{k}{4}^{\frac{2k+1}{8}}   e^{\sqrt{2}(X_{t+\ln 8}-X_{\ln 8})(x)-t}\,dx.
\end{align*}
We set
\begin{align*}
Y_i=&\frac{1}{8}\inf_{x\in [\frac{i}{4},\frac{2i+1}{8}]}e^{\sqrt{2}X_{\ln 8}(x)-\ln 8}\; , & N_i=8\int_\frac{i}{4}^{\frac{2i+1}{8}}   e^{\sqrt{2}(X_{t+\ln 8}-X_{\ln 8})(x)-t}\,dx.
\end{align*}
for $i=0,\dots,3$. A straightforward computation of covariances shows that $((X_{t+\ln 8}-X_{\ln 8})(x))_{x \in \R}$ has same distribution as $(X_{t}(8x))_{x \in \R}$. It is plain to deduce that:\\
$\bullet$  the random variables $(Y_i)_i$ are independent of $(N_i)_i$,\\
$\bullet$  the random variables $(Y_i)_i$ are identically distributed,\\
$\bullet$  the random variables $(N_i)_i$ are identically distributed with common law $M_t^{\sqrt{2}}([0,1])$.
 
\vspace{1mm} 
Let us define
$$\varphi_t(s)=\E[e^{-s\sqrt{t}M_t^{\sqrt{2}}([0,1])}].$$
Since the mapping $x\mapsto e^{-s\sqrt{t}x}$ is convex, we can apply Lemma \ref{lem:cvx} and obtain  $\varphi_t(s)\leq \E[e^{-s\sqrt{t}M_{t+\ln 8}^{\sqrt{2}}([0,1])}]$. Therefore:
\begin{align*}
\varphi_t(s)\leq   \E[e^{-s\sqrt{t}M_{t+\ln 8}^{\sqrt{2}}([0,1])}] 
\leq  \E[e^{-s\sqrt{t}\sum_{i=0}^3Y_iN_i}]=\E[\prod_{i=0}^3\varphi_t(sY_i) )]\leq \E[\varphi_t(sY_1)^4].
\end{align*}
Let us denote by $F$ the distribution function of $Y_1$. We have:
\begin{align*}
\varphi_t(s)\leq  &\int_0^\infty\varphi_t(sx)^4\,F(dx)=\int_0^{s^{-1/2}} \varphi_t(sx)^4\,F(dx)+\int_ {s^{-1/2}}^\infty \varphi_t(sx)^4\,F(dx)\\
\leq  &F(s^{-1/2})+ \varphi_t(s^{1/2})^4.
\end{align*}
 It is plain to check that $\E[(Y_1)^{-q}]<+\infty$ for all $q>0$: use for instance \cite[Lemma 5, section B]{bacry} to prove the result for the cone based process $\omega_l$ (with the notations of \cite{bacry}) and then extend the result for $X_t$ thanks to Lemma \ref{lem:cvx}.  Fix $q_0>0$. By Markov's inequality, we deduce:
\begin{equation*}
\varphi_t(s)\leq s^{-\frac{q_0}{2}}\E[(Y_1)^{-q_0}]+  \varphi_t(s^{1/2})^4 .
\end{equation*}
Therefore
\begin{equation}\label{est:phits}
\varphi_t(s)  \leq \big(s^{2q-\frac{q_0}{2}} \E[(Y_1)^{-q_0}]+  \varphi_t(s^{1/2})^2 \big) \big(s^{-2q} +  \varphi_t(s^{1/2})^2 \big).
\end{equation}
Let us admit for a while the following lemma
\begin{lemma}\label{cvunif}
The family of functions $(\varphi_t(\cdot))_t$  converges uniformly on $\R_+$  as $t\to \infty$ towards
$$\varphi(s)=\E[e^{-s\sqrt{2/\pi}M'([0,1])}].$$
\end{lemma}
Since $\P(M'([0,1])=0)=0$, we deduce $\lim_{s\to \infty}\varphi(s)=0$. Therefore, for  $q<q_0/4$, there exist $s_0>0$ and $t_0>0$ such that $$\forall s>s_0,\forall t>t_0 \quad \big(s^{2q-\frac{q_0}{2}} \E[(Y_1)^{-q_0}]+  \varphi_t(s^{1/2})^2 \big)\leq \frac{1}{\sqrt{2}}.$$
Plugging this estimate into \eqref{est:phits} yields:
\begin{equation*}
\forall s>s_0,\forall t>t_0 \quad \varphi_t(s)  \leq  \frac{1}{\sqrt{2}}\big(s^{-2q} +  \varphi_t(s^{1/2})^2 \big).
\end{equation*}
Then, by induction and by using $(a+b)^2\leq 2a^2+2b^2$, we check that for $s>s_0$:
\begin{equation}\label{ind}
\varphi_t(s^{2^n})\leq  \frac{1}{\sqrt{2}}\big(a_{n+1}s^{-q2^{n+1}}+\varphi_t^{2^{n+1}}(s^{1/2})\big)
\end{equation} where $a_1=1$ and $a_{n+1}=a_n^2+1$. Let us choose $Q\geq \frac{1+\sqrt{5}}{2}$ so as to have $Q^2- Q-1\geq 0$. It is then plain to check by induction that $Q^{2^n}\geq Q+a_n$, and thus $Q^{2^n}\geq a_n$. Then for all $x>s$, there exists $n\in\N$ such that $s^{2^n}\leq x<s^{2^{n+1}}$. In other words,
$$2^n\leq \frac{\ln x}{\ln s}<2^{n+1} .$$
Thus we obtain from \eqref{ind}
\begin{align*}
\varphi_t(x)\leq &\varphi(s^{2^n})\leq   \frac{1}{\sqrt{2}}\big(Q^{2^{n+1}}s^{-q2^{n+1}}+\varphi_t^{2^{n+1}}(s^{1/2})\big) \leq    \frac{1}{\sqrt{2}}(x^{-\alpha}+x^{-\beta})
\end{align*}
where $\alpha=q-\frac{\ln Q}{\ln s}>0$ and $\beta=-  \frac{\ln \varphi_t(s^{1/2})}{\ln s}>0$. Let us choose $s=s_0+Q^{1/q}$. Since $\lim_{t\to \infty}\varphi_t(s)=\E[e^{-sM'([0,1])}]$, we can choose $\beta$ arbitrarily close to $\beta_0=- \frac{\ln \varphi(s^{1/2})}{\ln s}>0$. To sum up, we have proved that, for all $\beta<\beta_0$, there exists $x_0=s_0+Q^{1/q}$ and $t_0>0$ such that
\begin{equation}\label{phit}
\forall x>x_0,\forall t>t_0 \quad \varphi_t(x)  \leq  \frac{1}{\sqrt{2}}(x^{-\alpha}+x^{-\beta}).
\end{equation}
It is straightforward to deduce
$$\sup_{t \geq 1}\E\Big[\Big(\frac{1}{\sqrt{t}M_t^{\sqrt{2}}(A)}\Big)^q\Big]<+\infty$$ for all $0<q<\min(\alpha,\beta_0)$. Put in other words, we have proved the result ``only for small $q$". But, remembering that
$$\varphi_t(x)\leq \E[\varphi_t(xY_1)^4]$$ and $\E[Y_1^{-q}]<+\infty$ for all $q>0$, we can deduce from \eqref{phit} the result for arbitrary $q$ by induction and by using Markov's inequality.\qed

\vspace{2mm}
\noindent {\it Proof of Lemma \ref{cvunif}.} Define the family of functions
$$f_t: y\in [0,1[\mapsto \varphi_t(\tan \frac{\pi y}{2}).$$ Since $\P\big(\sqrt{t}M_t^{\sqrt{2}}([0,1])=0\big)=0$, it is plain to deduce that $f_t$ can be continuously extended to $[0,1]$ by setting $f_t(1)=0$ for all $t$. In the same way we define a continuous function on $[0,1]$  by
$$f: y\in [0,1[\mapsto \varphi(\tan \frac{\pi y}{2})$$ and $f(1)=0$ (possible because $\P\big(M'([0,1])=0\big)=0$).  The family $(f_t)_t$ pointwise converges as $t\to \infty$ towards $f$. Furthermore, the functions $f_t$ are non increasing for all $t$. It is then standard to deduce the uniform convergence. The lemma follows. \qed

\vspace{3mm}
\noindent {\it Proof of Corollary \ref{prop:pw}}. Consider a solution $M$ of the $\star$-equation \eqref{star} with $\omega_\varepsilon $ given by \eqref{staromega} and $\gamma^2=2d$.
Consider an exponent $0<q<1$, an open bounded set $A$ and $\lambda<1$. From equation \eqref{star} (with $\varepsilon=\lambda$) and Jensen's inequality we have:
\begin{align*}
 \E[(M(\lambda A)^q]&= \E\Big[\Big( \int_{\lambda A}e^{\gamma X_{\ln \frac{1}{\lambda}}(r)-\frac{\gamma^2}{2}\E[X_{\ln \frac{1}{\lambda}}(r)^2]}\lambda^d M^ {\lambda}(dr)\Big)^q\Big]\\
 &\geq \lambda^{dq}\E\Big[e^{\gamma q X_{\ln \frac{1}{\lambda}}(0)-\frac{\gamma^2q}{2}\E[X_{\ln \frac{1}{\lambda}}(0)^2]}   \widetilde{M}(A) ^q\Big]\\
 &=\lambda^{\xi(q)}\E\big[   M (A) ^q\big]
\end{align*}
where $ \widetilde{M}$ is a copy of $M$ independent of $X$.
Conversely, we use Lemma \ref{lem:cvx} and equation \eqref{star}:
\begin{align*}
 \E[M(\lambda A)^q]&= \E\Big[\Big( \int_{\lambda A}e^{\gamma X_{\ln \frac{1}{\lambda}}(r)-\frac{\gamma^2}{2}\E[X_{\ln \frac{1}{\lambda}}(r)^2]}M^{ \lambda}(dr)\Big)^q\Big]\\
&\leq \E\Big[\Big( \int_{\lambda A}e^{\gamma \sqrt{a_\lambda} Z-\frac{\gamma^2}{2}a_\lambda}M^{ \lambda}(dr)\Big)^q\Big]
\end{align*}
where $Z$ is a standard Gaussian random variable independent of $M^\lambda$ and $a_\lambda=\inf_{x,y\in\lambda A}K_{  \ln \frac{1}{\lambda}}(y-x)$. Therefore:
\begin{align*}
 \E[M(\lambda A)^q] &\leq \lambda^{dq}\E\big[e^{\gamma q \sqrt{a_\lambda} Z-\frac{q\gamma^2}{2}a_\lambda }\big]\E\Big[   M (A) ^q\Big] .
\end{align*}
We complete the proof by noticing that $a_\lambda\geq C+\ln\frac{1}{\lambda}$, which follows from the fact that $k$ is Lipschitz at $0$.\qed

\begin{remark}
Actually, using the $\star$-equation to compute the power-law spectrum is not necessary: it  can be computed  with similar arguments for any derivative multiplicative chaos associated to a log-correlated Gaussian field.
\end{remark}

\section{Proof of Section \ref{sec:KPZ}}

\noindent {\it Proof of Theorem \ref{KPZ}.}   Without loss of generality, we assume that $k$ vanishes outside the ball $B(0,1)$. Let $K$ be a compact set included in the ball $B(0,1)$ with Lebesgue Hausdorff dimension $ {\rm dim}_{Leb}(K)$.

We first want to prove $d\,{\rm dim}_{Leb}(K)\geq \xi({\rm dim}_M(K))$. If ${\rm dim}_{Leb}(K)=1$, the inequality is trivial because $\xi(q)\leq d$ for all $q$. So we may assume ${\rm dim}_{Leb}(K)<1$. Let $q\in [0,1[$ be such that $\xi(q)>d\,{\rm dim}_{Leb}(K)$ with $\xi(q)<d$. For $\varepsilon>0$, there is a covering of $K$ by a countable family of balls $(B(x_n,r_n))_n$ such that
$$ \sum_n r_n^{\xi(q)}<\varepsilon.$$ By using in turn  the stationarity and the power law spectrum of the measure, we have
\begin{align*}
\E\Big[\sum_nM'(B(x_n,r_n))^q\Big]& =\sum_n\E\Big[M'(B(0,r_n))^q\Big]\leq C_q \sum_n r_n^{\xi(q)}\leq C_q\varepsilon.
\end{align*}
Using Markov's inequality, we deduce
$$\P\Big(\sum_nM'(B(x_n,r_n))^q\leq C_q \sqrt{\varepsilon}\Big)\geq 1-\sqrt{\varepsilon}.$$
Thus, with probability $1-\sqrt{\varepsilon}$, there is a covering of balls of $K$ such that
$$\sum_nM'(B(x_n,r_n))^q\leq C_q\sqrt{\varepsilon}.$$ So $q\geq {\rm dim}_{M'}(K)$ almost surely. Therefore $d\,{\rm dim}_{Leb}(K)\geq \xi({\rm dim}_M(K))$.

Conversely, consider $q\in[0,1]$ such that $\xi(q)<d\,{\rm dim}_{Leb}(K)$. Observe that necessarily $q<1$ because $\xi(1)=d$ and ${\rm dim}_{Leb}(K)\leq 1$. So we may assume $q<1$. 
By  Frostman's Lemma, there is a probability measure $\gamma$ supported by $K$ such that
\begin{equation}\label{energy}
\int_{B(0,1)^2}\frac{1}{|x-y|^{\xi(q)}}\gamma(dx)\gamma(dy)<+\infty.
\end{equation}
For $q\in [0,1[$, let us define the random measure $\widetilde{\gamma}$ as the almost sure limit of the following family of positive random measures:
\begin{equation}\label{gamtilde}
\widetilde{\gamma}(dx)=\lim_{t \to \infty}e^{q\sqrt{2d} X_t(x) - q^2d\E[(X_t(x))^2]}\gamma(dx).
\end{equation}
The limit is non trivial: Kahane \cite{cf:Kah} (see also  a more general proof in \cite{sohier}) proved that, for a Radon measure $\gamma(dx)$ satisfying \eqref{energy} with a power exponent $\kappa$ (instead of $\xi(q)$ in \eqref{energy}) the associated chaos
$$\widetilde{\gamma}(dx)=\lim_{t\to \infty}e^{s X_t(x) - \frac{s^2}{2}\E[(X_t(x))^2]}\gamma(dx) $$ is non degenerate (i.e. the martingale is regular) provided that $\kappa-\frac{s^2}{2}>0$. In our context, this condition reads $q^2d<\xi(q)$, that is $q<1$.
%

 From   Frostman's lemma again (this form is non standard as it involves measures instead of distances, see \cite[Lemma 17]{BJRV} for a proof), we just have to prove that the quantity
\begin{equation}\label{frofro}
\int_{B(0,1)^2}\frac{1}{M'(B(x,|y-x|))^{q}}\widetilde{\gamma}(dx)\widetilde{\gamma}(dy)
\end{equation}
is finite almost surely. Actually \cite[Lemma 17]{BJRV} requires to consider each of the $2^d$ portions of the ball: $B(x,|y-x|)\cap \big\{z\in\R^d; (z_1-x_1)\varepsilon_1\geq 0,\dots,(z_d-x_d)\geq \varepsilon_d\geq 0\big\}$ for $\varepsilon_1,\dots,\varepsilon_d\in \{-1,1\}$. Yet, this is exactly the same proof but considering these portions of balls would heavily decrease the readability of the proof because of notational issues. So we keep on considering a ball in   \eqref{frofro}. Finiteness of \eqref{frofro} follows  from
\begin{equation}\label{froststar}
\E\Big[\int_{B(0,1)^2}\frac{1}{M'(B(x,|y-x|))^{q}}\widetilde{\gamma}(dx)\widetilde{\gamma}(dy) \Big]<+\infty.
\end{equation}
Actually, by using  Fatou's lemma in \eqref{froststar}, it is even sufficient to prove that:
\begin{equation}\label{frostt}
\liminf_{t\to \infty}  \int_{B(0,1)^2}\E\Big[\frac{e^{q\sqrt{2d} X_t(x)+q\sqrt{2d} X_t(y)- 2q^2d \E[(X_t(x))^2]}}{\big(\sqrt{t}M^{\sqrt{2d}}_t(B(x,|y-x|))\big)^{q}} \Big] \gamma(dx) \gamma (dy)<+\infty.
 \end{equation}
From now on, we will focus on computing the above integral \eqref{frostt}. There is a way of making the computations with minimal effort: we exchange the process $X_t$ with the perfect scaling process introduced in \cite{rhovar}. We just have to justify that this change of processes is mathematically rigorous. So let us admit for a while the following lemma:

\begin{lemma}\label{changeofproc}
If \eqref{frostt} is finite for the process $(X_t(x))_{t,\geq 0,x\in\R^d}$ with correlations given by
\begin{equation}\label{newcorr}
\E[X_t(x),X_s(y)]=\int_Sg_{e^{-\min(s,t)}}(|\langle x-y,s\rangle|)\sigma(ds)
 \end{equation}
 where $S$ stands for the sphere of $\R^d$, $\sigma$ the uniform measure on the sphere and the function $g_u$ (for $0<u\leq 1$) is given by
$$g_u(r)=\left\{\begin{array}{ll}
\ln_+\frac{2}{r} & \text{if }r\geq u\\
\ln \frac{2}{u}+1-\frac{r}{u} & \text{if }r< u
\end{array}\right. $$
then \eqref{frostt} is finite for the process $(X_t(x))_{t,\geq 0,x\in\R^d}$ with correlations given by \eqref{corrX}.
\end{lemma}

So, from now on, we assume that the correlations of $(X_t(x))_{t,\geq 0,x\in\R^d}$ are the new correlations specified in Lemma \ref{changeofproc} (see \cite{rhovar} for further details). Notice also that the measure $M^\gamma$ also involves this new process. Such a family of kernels possesses useful scaling properties, namely that for $|x|\leq 1$ and $h,t\geq 0$, $K_{t+h }(e^{-h} x)=K_t(x)+h$. In particular, we have the following scaling relation for  all   $h,t\geq 0$:
\begin{align}\label{scalingkl}
 \big((X_{t+h}(e^{-h} x))_{x\in B(0,1)}&,(M^{\sqrt{2d}}_{t+h}(e^{-h} A))_{A\subset B(0,1)}\big)\nonumber\\
 &\stackrel{law}{=}\big((X_t(x)+\Omega_h)_{x\in B(0,1)},( e^{\sqrt{2d}\Omega_h-2dh}M^{\sqrt{2d}}_{  t}(  A))_{A\subset B(0,1)}\big).
\end{align}
where $\Omega_h$ is a centered Gaussian random variable with variance $h$ and independent of the couple $$\big((X^t(x))_{x\in B(0,1)},( M^{\sqrt{2d}}_{  t}(  A))_{A\subset B(0,1)}\big).$$ We will use the above relation throughout the proof.

By stationarity of the process $X_t$,  we can translate the integrand in the quantity \eqref{frostt} to get
\begin{align*}
\int_{B(0,1)^2}&\E\Big[\frac{e^{q\sqrt{2d} X_t(x)+q\sqrt{2d} X_t(y)- 2q^2d \E[(X_t(x))^2]}}{\big(\sqrt{t}M^{\sqrt{2d}}_t(B(x,|y-x|))\big)^{q}} \Big] \gamma(dx) \gamma (dy) \\
&\leq     \int_{B(0,1)^2}\E\Big[\frac{e^{q\sqrt{2d} X_t(0)+q\sqrt{2d} X_t(y-x)- q^2d \E[(X_t(0))^2]}}{\big(\sqrt{t}M^{\sqrt{2d}}_t(B(0,|y-x|))\big)^{q}}\Big] \gamma(dx) \gamma (dy).
 \end{align*}
Now we split the latter integral according to scales larger or smaller than $e^{-t}$ and obtain
\begin{align*}
\int_{B(0,1)^2}\E\Big[&\frac{e^{q\sqrt{2d}  X_t(0)+q\sqrt{2d} X_t(y-x)- q^22d \E[(X_t(0))^2]}}{\big(\sqrt{t}M^{\sqrt{2d}}_t(B(0,|y-x|))\big)^{q}}\Big] \gamma(dx) \gamma (dy)\\
=  & \int_{|x-y|>e^{-t}}\E\Big[\frac{e^{q\sqrt{2d} X_t(0)+q\sqrt{2d} X_t(y-x)- q^22d \E[(X_t(0))^2]}}{\big(\sqrt{t}M^{\sqrt{2d}}_t(B(0,|y-x|))\big)^{q}}\Big] \gamma(dx) \gamma (dy)\\
&+ \int_{|x-y|\leq e^{-t}}\E\Big[\frac{e^{q\sqrt{2d} X_t(0)+q\sqrt{2d} X_t(y-x)- q^22d \E[(X_t(0))^2]}}{\big(\sqrt{t}M^{\sqrt{2d}}_t(B(0,|y-x|))\big)^{q}}\Big] \gamma(dx) \gamma (dy)\\
\stackrel{def}{=}&I_1(t)+I_2(t)
 \end{align*}

We first estimate $I_1(t)$.
 For $|x-y|>e^{-t}$, by using \eqref{scalingkl} with $e^{-h}=|y-x|$, we deduce
\begin{align*}
\E\Big[&\frac{e^{q\sqrt{2d} X_t(0)+q\sqrt{2d} X_t(y-x)- q^22d \E[(X_t(0))^2]}}{\big(\sqrt{t}M^{\sqrt{2d}}_t(B(0,|y-x|))\big)^{q}}\Big] \\
&=\E\Big[ \frac{e^{q2\sqrt{2d} \Omega_h - 2q^2d h}e^{q\sqrt{2d} X_{t-h}(0)+q\sqrt{2d} X_{t-h}( \frac{y-x}{|y-x|} )- 2q^2d \E[(X_{t-h}(0))^2]}}{\Big(  e^{\sqrt{2d}\Omega_h-2dh}\sqrt{t} M^{\sqrt{2d}}_{  t-h}(B(0,1))     \Big)^q }\Big]  \\
&= \E\Big[ e^{q \sqrt{2d} \Omega_h - 2(q^2-q)d h}\Big]\E\Big[\frac{e^{q\sqrt{2d} X_{t-h}(0)+q\sqrt{2d} X_{t-h}( \frac{y-x}{|y-x|} )- 2q^2d \E[(X_{t-h}(0))^2]}}{\Big(  \sqrt{t} M^{\sqrt{2d}}_{  t-h}(B(0,1))     \Big)^q }\Big]  \\
&\leq    \frac{C}{|y-x|^{\xi(q)}}\E\Big[\frac{ 1}{\Big(  \sqrt{t} M^{\sqrt{2d}}_{  t-h}(B(0,1/2))      \Big)^q } \Big].
  \end{align*}
In the last line, we have used a Girsanov transform to get rid of the numerator. Therefore, if we set $c=\sup_{u \leq 1}\E\Big[\frac{ 1}{\Big( M^{\sqrt{2d}}_{  u}(B(0,1/2))      \Big)^q } \Big]$, then
\begin{equation}\label{estI1}
I_1(t)\leq C( \sup_{u \geq 1}\E\Big[\frac{ 1}{\Big(  \sqrt{u} M^{\sqrt{2d}}_{  u}(B(0,1/2))      \Big)^q } \Big]+\frac{c}{\sqrt{t}}) \int_{|x-y|>e^{-t}}\frac{C}{|y-x|^{\xi(q)}}\gamma(dx)\gamma(dy).
\end{equation}
From Corollary \ref{momneg}, the quantity $$\sup_{u \geq 1}\E\Big[\frac{ 1}{\Big(  \sqrt{u} M^{\sqrt{2d}}_{  u }(B(0,1/2))      \Big)^q } \Big]  $$ is finite.

\begin{remark}
Corollary \ref{momneg} only deals with Gaussian fields with correlations given by \eqref{corrX}. From Kahane's convexity inequality, this quantity is also finite for  Gaussian fields with correlations given by $$K(x,y)= \ln_+\frac{1}{|y-x|}+g(x,y)$$ for some bounded function $g$. In particular, it is finite for the field considered in Lemma \ref{changeofproc}.
\end{remark}

To treat the term $I_2(t)$, we use quite a similar argument except that we use the scaling relation on $h=t$ instead of $-\ln |y-x|$, and the Girsanov transform again:
\begin{align*}
 &I_2(t)\\&=    \int_{|y-x|\leq e^{-t}}\E\Big[\frac{e^{2q\sqrt{2d} \Omega_{t}-q^22dt}e^{q\sqrt{2d}  X_{0}(0)+q\sqrt{2d}  X_0(e^t(y-x)) -q^22d\E[(X_0(0))^2]}}{\Big(e^{q\sqrt{2d}\Omega_{t}-q2dt }\sqrt{t}M_{0}\big(B(0,e^t|y-x|)\big)\Big)^{q}}\Big] \gamma(dx) \gamma (dy)\\
&=  \int_{|y-x|\leq e^{-t}}\E\Big[ e^{  q\sqrt{2d} \Omega_{t}-  2(q^2-q)d t}\Big]\E\Big[\frac{e^{q\sqrt{2d}  X_{0}(0)+q\sqrt{2d}  X_0(e^t(y-x)) -q^22d\E[(X_0(0))^2]}}{ \Big( \sqrt{t}M_{0}\big(B(0,e^t|y-x|)\big)\Big)^{q}}\Big] \gamma(dx) \gamma (dy)\\
&=    \int_{|y-x|\!\leq\! e^{-t}}\!\!\frac{e^{t\xi(q)}}{ t^{q/2}}\E\Big[\frac{e^{q^22dK_0\big(e^t(y-x)\big) }}{ \Big(\!\!\int_{B(0,e^t|y-x|)}
e^{ \sqrt{2d}X_{0}(u)-d\E[(X_0(0))^2]+q2dK_0(e^t(y-x)-u)+q2dK_0(u)}\,du\!\Big)^q}\Big] \gamma(dx) \gamma (dy).
 \end{align*}
By using the fact that $K_0$ is positive and bounded by $\ln 2$, we have (for some positive constant $C$ independent of $t$)
$$I_2(t)\leq C\int_{|y-x|\leq e^{-t}}\frac{e^{t\xi(q)}}{ t^{q/2}} \E\Big[\frac{1}{ \Big(\int_{B(0,e^t|y-x|)}
e^{ \sqrt{2d}X_{0}(u)-d\E[(X_1(0))^2]}\,du\Big)^q}\Big] \gamma(dx) \gamma (dy).$$
Since $\E[X_0(u)X_0(0)]\leq \E[(X_0(0))^2]$, we can use Kahane's convexity inequalities to the convex mapping $x\mapsto \frac{1}{x^q}$. We deduce (for some positive constant $C$ independent of $t$, which may change from line to line)
\begin{align*}
I_2(t)\leq &C\int_{|y-x|\leq e^{-t}}e^{t\xi(q)} t^{-q/2}\E\Big[\frac{1}{ \Big(\int_{B(0,e^t|y-x|)}
e^{ \sqrt{2d}X_{0}(0)-d\E[(X_1(0))^2]}\,du\Big)^q}\Big]\gamma(dx) \gamma (dy)\\
\leq & C t^{-q/2} \int_{|y-x|\leq e^{-t}}e^{t\xi(q)}           \frac{1}{e^{dqt}|y-x|^{dq} }  \gamma(dx) \gamma (dy)\\
\leq &C t^{-q/2}    \int_{|y-x|\leq e^{-t}}\frac{1}{|y-x|^{\xi(q)} }  \gamma(dx) \gamma (dy).
\end{align*}
Hence $$\lim_{t \to \infty}I_2(t)\leq C\limsup_{t\to\infty}\,\,\,t^{-q/2}\int_{|y-x|\leq e^{-t}}\frac{1}{|y-x|^{\xi(q)} }  \gamma(dx) \gamma (dy)=0.$$
The KPZ formula is proved.\qed

\vspace{2mm}
\noindent {\it Proof of Lemma \ref{changeofproc}.} Let us denote by $K_t$ the $\star$-scale invariant kernel given by \eqref{corrX} associated to the process $(X_t(x))_{t,x}$. We will use the superscript $^p$ to denote the corresponding quantities associated to the ``perfect" kernel of \cite{rhovar}: we denote by $K_t^p$ the kernel described in Lemma \ref{changeofproc}, by $(X_t^p(x))_{t,x}$ (resp. $M^{p,\sqrt{2d}}_t$) the associated Gaussian field (resp. approximate multiplicative chaos). It is plain to see there is a constant $C>0$ such that
\begin{equation}\label{relperf}
\forall t>0, \forall x\in\R^d,\quad K_t(x)-C\leq K_t^p(x)\leq K_t(x)+C.
\end{equation}

Now we prove the Lemma. By using the Girsanov transform, we have:
\begin{align*}
& \int_{B(0,1)^2}\E \Big[\frac{e^{q\sqrt{2d} X_t(x)+q\sqrt{2d} X_t(y)- 2q^2d \E[(X_t(x))^2]}}{\big(\sqrt{t}M^{\sqrt{2d}}_t(B(x,|y-x|))\big)^{q}} \Big] \gamma(dx) \gamma (dy)\\
 &=\int_{B(0,1)^2}\E\Big[\frac{1}{\big(\sqrt{t}\int_{B(x,|y-x|)}e^{q2dK_t(u-x)+q2dK_t(u-y)}M^{\sqrt{2d}}_t(du)\big)^{q}} \Big]e^{q^22dK_t(x-y)} \gamma(dx) \gamma (dy)\\
&\leq  e^{(2d q^2+4dq)C}\int_{B(0,1)^2}\E\Big[\frac{1}{\big(\sqrt{t}\int_{B(x,|y-x|)}e^{q2dK^p_t(u-x)+q2dK^p_t(u-y)}M^{\sqrt{2d}}_t(du)\big)^{q}} \Big]e^{q^22dK_t^p(x-y)} \gamma(dx) \gamma (dy).
 \end{align*}
 By using Kahane's convexity inequality (Lemma \ref{lem:cvx}) together with the inequality \eqref{relperf} between covariance kernels, we get:
\begin{align*}
\int_{B(0,1)^2}&\E\Big[\frac{1}{\big(\sqrt{t}\int_{B(x,|y-x|)}e^{q2dK^p_t(u-x)+q2dK^p_t(u-y)}M^{\sqrt{2d}}_t(du)\big)^{q}} \Big]e^{q^22dK_t^p(x-y)} \gamma(dx) \gamma (dy)\\
&\leq C'\int_{B(0,1)^2} \E\Big[\frac{1}{\big(\sqrt{t}\int_{B(x,|y-x|)}e^{q2dK^p_t(u-x)+q2dK^p_t(u-y)}M^{p,\sqrt{2d}}_t(du)\big)^{q}} \Big]e^{q^22dK_t^p(x-y)} \gamma(dx) \gamma (dy).
 \end{align*}
Then, by  the Girsanov transform again, we have
\begin{align*}
\int_{B(0,1)^2}\E&\Big[\frac{1}{\big(\sqrt{t}\int_{B(x,|y-x|)}e^{q2dK^p_t(u-x)+q2dK^p_t(u-y)}M^{p,\sqrt{2d}}_t(du)\big)^{q}} \Big]e^{q^22dK_t^p(x-y)} \gamma(dx) \gamma (dy)\\
&=\int_{B(0,1)^2}\E \Big[\frac{e^{q\sqrt{2d} X^p_t(x)+q\sqrt{2d} X^p_t(y)- 2q^2d \E[(X^p_t(x))^2]}}{\big(\sqrt{t}M^{p,\sqrt{2d}}_t(B(x,|y-x|))\big)^{q}} \Big] \gamma(dx) \gamma (dy).
 \end{align*}
The lemma follows.\qed

\vspace{1mm}
To sum up, we have proved that establishing the KPZ formula for the perfect kernel is equivalent to establishing the KPZ formula for all $\star$-scale invariant kernels. Furthermore the above argument is obviously valid for any log-correlated Gaussian field (as it only involves the Girsanov transform) and  for other values of $\gamma$: for $\gamma^2< 2d$ with the techniques developed in \cite{cf:RhoVar} or for $\gamma^2>2d$ with the techniques developed in \cite{BJRV}. In particular and as advocated in \cite{BJRV,cf:RhoVar}, the KPZ formula established in \cite{BJRV,cf:RhoVar} in terms of {\bf Hausdorff dimensions} are valid for the GFF and also for all log-correlated Gaussian fields in any dimension.  The reader may compare with \cite{cf:DuSh} where the KPZ formula is stated in terms of {\bf expected box counting dimensions}. At criticality, things are a bit more subtle: the KPZ formula established in this paper is valid for all the derivative multiplicative chaos for which you can establish the renormalization theorem \ref{seneta}. This theorem is necessary to be in position to apply Kahane's convexity inequalities. In particular,
the above proof also establishes Theorem \ref{KPZGFF} (and for the MFF as well).

\section{Proofs for Free Fields}\label{freefields}
In this section, we  extend our results to the case of free fields.

\subsection{Whole-Plane Massive Free Field}
We consider a white noise decomposition $(X_t)_t$ with either of the following covariance structures:
\begin{align}
\label{covMFF1}\E[X_t(x)X_s(y)]=&\int_1^{e^{\min(s,t)}}\frac{k_m(u(x-y))}{u}\,du:=G_{\min(s,t)}(x-y),\\
\label{covMFF2}\E[X_t(x)X_s(y)]=&\pi \int^{\infty}_{e^{-2\min(s,t)}} p(u,x,y)e^{-\frac{m^2}{2}u}\,du:=G_{\min(s,t)}(x-y),
\end{align}
with $p(u,x,y)=\frac{e^{-\frac{|x-y|^2}{2 u}}}{2 \pi u} $ and $k_m(z)=\frac{1}{2}\int_0^\infty e^{-\frac{m^2}{2v}|z|^2-\frac{v}{2}}\,dv$.  
We will explain below how to construct the derivative multiplicative chaos and how to prove the Seneta-Heyde renormalization for both of them.

\vspace{2mm}
{\bf Derivative multiplicative chaos.} Concerning the construction of the derivative multiplicative chaos, we can go along the lines of \cite{DRSV}. Yet, some observations help to manage the long-range correlations. Sticking to the notations of \cite{DRSV}, the point is to prove that the martingale $(Z^\beta_t(A))_t$ is regular (and then to prove that the limiting measure has no atom), i.e., to follow the proof of \cite[Proposition 14]{DRSV}. Basically, all computations rely on the joint law of the couple of processes $(X_t(w),X_t(x))_{t\geq 0} $ for $x\not =w$: one has two quantities to estimate, $\widetilde{\Pi}_1$ and $\widetilde{\Pi}_2$, and one can see that the following computations only rely on the law of this couple in equations (19) and (29) in \cite[proof of Prop.14]{DRSV}. (These conditional expectations are then integrated with respect to the law of $(X_t(x))_t$, so that considering their law indeed suffices.) The key point in \cite{DRSV} is the following: if the kernel $k$ has compact support included in the ball $B(0,1)$ then  the fields  $\big(X_t(w)-X_{s_0}(w),X_t(x)-X_{s_0}(x)\big)_{t\geq s_0} $ are independent for $s_0\geq\ln\frac{1}{|x-w|}$. 

In the case of the MFF, this decorrelation property clearly does not hold. The idea is then to encode the long-range correlations into a nice field that can be easily controlled, and the logarithmic short-range correlations into another field possessing the decorrelation property. Let us describe the case where the correlations are given by \eqref{covMFF1}. We introduce $4$ independent Brownian motions $B^i$ for $i=1,\dots,4$, and a smooth function $\varphi:\R_+\to\R$ such that $0\leq \varphi\leq 1$, $\varphi=1$ over $[0,1/2]$ and $\varphi=0$ over $[2,+\infty[$. We denote by $G'_t(x)$ the derivative at time $t$  of the mapping $t\mapsto G_t(x)$. (This derivative is positive.)  Then we  define
\begin{align*}
P^{sc}_t=&\int_0^t \sqrt{\varphi(e^r|x-w|)G'_r(x-w)}\,dB_r^1,\quad &P^{lc}_t=&\int_0^t\sqrt{(1-\varphi(e^r|x-w|))G'_r(x-w)}\,dB^2_r, \\
Z^{x}_t=&\int_0^t\sqrt{G'_r(0)-G'_r(x-w)}\,dB^3_r,\quad &Z^{w}_t=&\int_0^t\sqrt{G'_r(0)-G'_r(x-w)}\,dB^4_r,
\end{align*}
and we have the following relation in law:
$$(X_t(w),X_t(x))_{t\geq 0}\stackrel{law}{=}\Big(P^{sc}_t+ P^{lc}_t+Z^w_t,P^{sc}_t+ P^{lc}_t+Z^x_t\Big)_{t\geq 0}.$$
The Gaussian process $P^{lc}$ encodes the long-range correlations and has independent temporal increments. Furthermore
\begin{align*}
\E[(P^{lc}_t)^2] =&\int_0^t(1-\varphi(e^r|x-w|))G'_r(x-w)\,dr= \int_0^t(1-\varphi(e^r|x-w|))k_m(e^r(x-w))\,dr\\ \leq & \int_0^\infty(1-\varphi(u))\frac{k_m(u \frac{x-w}{|x-w|})}{u} \,du.
\end{align*}
Hence, we deduce that $\sup_{x\not =w\in \R^2,t\geq 0}\E[(P^{lc}_t)^2]<+\infty$, therefore there exists a constant $C>0$ independent of $x,y$, such that
$\sup_{t \geq 0} (-P^{lc}_t)$ is stochastically dominated by $\sup_{0 \leq t  \leq C} B_t$, where $B_t$ is  standard Brownian motion. Hence there exists $\lambda,C>0$ such that
$$ \sup_{x\not =w\in \R^2}\E[e^{\lambda (\sup_{t\geq 0}(-P^{lc}_t))^2}] \leq C.$$
Therefore, one may easily get rid of this process in the proof of \cite[Proposition 14]{DRSV}: just condition on the process $\sup_{t\geq 0}(-P^{lc}_t)$. Then one just has to deal with the couple $\Big(P^{sc}_t+Z^w_t,P^{sc}_t+Z^x_t\Big)_{t\geq 0}$, whose increments are decorrelated from $s_0=\ln\frac{1}{|x-w|}$ on. One can then follow the rest of the proof in \cite[Proposition 14]{DRSV} up to harmless modifications. 

This argument may be adapted to quite a general family of log-correlated fields. For instance, in the case of correlations given by \eqref{covMFF2},  apply the same strategy, except for replacing now $\varphi(e^r|x-w|)$ by $\varphi(r^{-1/2}|x-w|) $: this different time parameterization just comes from the fact that the time evolution of the function $G_t$ in \eqref{covMFF2} differs from that in \eqref{covMFF1}.\qed

\vspace{2mm} 
{\bf Extension of the Seneta-Heyde renormalization.} The Seneta-Heyde renormalization requires another argument. We will use coupling arguments for the MFF and GFF. First we extend our Theorem \ref{seneta} to the following set-up:
\begin{theorem}\label{senetapower}
 Assume that the family $(X_t(x))_{t\geq 0,x\in\R^d}$ has covariance kernel of the form
\begin{equation}\label{corralpha}
\E[X_t(x)X_s(y)]=K_{\min(s,t)}(x-y),\quad \text{with }\quad K_t(x)=\int_1^{e^t}\frac{k(v,x)}{v}\,dv
\end{equation}
for some family $(k(v,\cdot))_{v\geq 1}$ of stationary covariance kernels such that:
 \begin{description}
\item[B.1] $k(\cdot,\cdot)$ is nonnegative,  continuous and $\int_1^{\infty}\frac{|1-k(v,0)|}{v}\,dv<+\infty$,
\item[B.2] for some constant $C$, for all $x\in\R^2$ and $v\geq 1$: $|k(v,0)-k(v,x)|\leq C|x|v$,
\item[B.3] $k(v,x)=0$ if $|x|\geq \psi(v)$, where $$\psi(v)=2\sqrt{\phi(v^{-2})}=Dv^{-1}(1+2\ln v)^\alpha$$ for some $\alpha>0$ and some fixed constant $D$.
\item[B.4] there exists $f:[1,+\infty[\to\R_+$ such that $k(v,x)\leq f(v|x|)$ for all $v\geq 1$ and $x\in\R^2$.
\end{description}
Then the conclusions of Theorem \ref{seneta} hold. 

\end{theorem} 
 
Let us admit for a while the validity of this theorem and make some comments. The kernel $k(v,x)$ plays the former role of $k(vx)$: so  basically,  Theorem \ref{seneta} corresponds to the case $\alpha=0$ and the above theorem slightly enlarges the allowed  size of the support of $k(v,x)$. This dependence of the size of the support does not allow us to deal with a covariance kernel $K(x)=\lim_{t\to\infty}K_t(x)$ with long-range correlations. The subtle thing here is that  a large class of kernels with long-range correlations may be decomposed into a  part satisfying Theorem \ref{senetapower}, plus a smooth (non logarithmically divergent) Gaussian field encoding the long-range correlations. So, in a way, the extra leeway corresponding to the additional logarithmic factor in [B.3] is what is needed to treat long-range correlations. We will apply this argument to the whole-plane MFF with the two possible types of  correlation structures, \eqref{covMFF1} or \eqref{covMFF2}. Observe that considering a well-chosen family of Gaussian processes with covariance structure \eqref{covMFF1} or \eqref{covMFF2} is not a restriction. Indeed, if we prove for one family $(X_t)_t$ the limit $\frac{\sqrt{t}M^{\sqrt{2d}}_t(A)}{M'_t(A)}\to \sqrt{2/\pi}$ in probability as $t\to+\infty$, then for any other family $(X_t)_t$ we will have the same convergence, at least in law, and since   $\sqrt{2/\pi}$ is a constant, we will also get convergence in probability for this family. So in what follows, we can construct as we please a family $(X_t)_t$ of Gaussian processes with covariance \eqref{covMFF1} or \eqref{covMFF2}.
\qed

\vspace{2mm} 
{\bf Application to the cut-off family with covariance \eqref{covMFF1}.} 

Let us set $G_m(x)=\lim_{t\to\infty}G_t(x)$. First, observe that 
\begin{align*}
G_m\star G_m(x)=&\int_{\R^2}G_m(y)G_m(y-x)\,dy= \pi^2\int_{\R^2}\int_0^{\infty}\int_0^{\infty}e^{-\frac{m^2(u+v)}{2}}p(u,0,y)p(v,x,y)\,dudvdy\\
=& \int_0^{\infty}\int_0^{\infty}e^{-\frac{m^2(u+v)}{2}}p(u+v,0,x) \,dudv=\pi^2\int_0^\infty vp(v,0,x)e^{-\frac{m^2v}{2}}\,dv\\
=&\frac{\pi}{m^2}k_m(x).
\end{align*}
Hence $G_m$ is a square root convolution of $k_m$, up to a multiplicative constant. Let us set $g_m(x)=\frac{m}{\sqrt{\pi}}G_m(x)$ in such a way $g_m\star g_m(x)=k_m(x)$. Let us set for $\alpha>1$
$$a_v:=(1+2\ln v)^{-\alpha},\quad \psi(v):=4\frac{(1+2\ln v)^\alpha}{v}= \frac{4}{v a_v}.$$ Consider a smooth function $\varphi$ defined on $\R$ such that $0\leq \varphi\leq 1$, $\varphi=1$ over $[0,1]$ and $\varphi=0$ on $[2,+\infty[$, no matter its values on $\R_-$. Then,  consider a white noise $W$ on $\R_+ \times \R^2$ and  define
\begin{align*}
X^{sc}_t(x)=&\int_1^{e^t} \int_{\R^2} \frac{g_m(y-xv)}{\sqrt{v}}\varphi(a_v |y-xv|)W(dv,dy)\\
X^{lc}_t(x)=&\int_1^{e^t} \int_{\R^2} \frac{g_m(y-xv)}{\sqrt{v}}(1-\varphi) (a_v |y-xv| )W(dv,dy),
\end{align*}
in such a way that the family $(X^{sc}_t+X^{lc})_{t\geq 0}$ has the same law as $(X_t)_{t\geq 0}$. The family $(X^{sc}_t)_t$ encodes the short-range  logarithmic divergence, whereas the family $(X^{lc}_t)_t$ encodes the long-range correlations.  Let us compute the covariance of $X^{sc}_t$:
\begin{align*}
K^{sc}_t(x):=\E[X^{sc}_t(0)X^{sc}_t(x)]=\int_1^{e^t}\int_{\R^2}g_m(y-xv)g_m(y)\varphi (a_v |y-xv| )\varphi (a_v |y| )dy\,\frac{dv}{v},
\end{align*}
in order to check that it satisfies the assumptions of Theorem \ref{senetapower}. Thus we set 
$$k(v,x):=\int_{\R^2} g_m(y-xv)g_m(y)\varphi (a_v |y-xv| )\varphi (a_v |y| )\,dy$$ in order to have $K^{sc}_t(x)=
\int_1^{e^t}k(v,x)\,\frac{dv}{v}$. First observe that
$$k(v,x)\leq \int_{\R^2} g_m(y-xv)g_m(y)\,dy=k_m(xv),$$ so that we can take $f(y)=k_m(y)$ to check [B.4]. Furthermore, it is obvious to see that $k(v,x)=0$ for $|x|\geq \psi(v)$, hence [B.3]. Then:
\begin{align*}
|1-k(v,0)|=&\int_{\R^2}g_m^2(y)(1-\varphi^2(a_v |y| )) \,dy 
\leq  \int_{\R^2\setminus B(0,(1+2\ln v)^\alpha)} g_m^2(y)\,dy. 
\end{align*}
Now, we use the fact that for some constant $C$ we have $g_m(y)\leq Ce^{-|y|}$ for $|y|\geq 2$, to deduce that [B.1] holds. To prove [B.2], we show a stronger statement: the uniform convergence as $v\to\infty$ of $k(v,\cdot/v)$ as well as its gradient towards $k_m$ and $\partial k_m$, with a control of the rate of convergence. As a byproduct, we will also get a uniform Kolmogorov criterion for $(X^{lc}_t)$. Consider $x\in\R^2$, we have
$$|k_m(x)-k(v,x/v)|= \int_{\R^2}g_m (y)g_m(x-y)\big(1-\varphi (a_v |y| )\varphi (a_v |y-x| )\big)\,dy .$$
If $|x|>a_v^{-1}/4$, a rough estimate (namely $1-\varphi\leq 1$) shows that the above quantity is less than $k_m(x)\leq k_m(a_v^{-1}/4)$. Otherwise $|x|\leq a_v^{-1}/4$ and in this case the product $\varphi (a_v |y| )\varphi (a_v |y-x| )$ is worth $1$ on the ball $B(0,a_v^{-1}/2)$, in which case we estimate the difference $|k_m(x)-k(v,x/v)|$ by 
\begin{align*}
|k_m(x)-k(v,x/v)|\leq & \int_{\R^2\setminus B(0,a_v^{-1}/2)}g_m (y)g_m(x-y)\,dy\\
\leq &C^2\int_{\R^2\setminus B(0,a_v^{-1}/2)}e^{-|y|}e^{-|x-y|}\,dy\\
\leq & C^2e^{-a_v^{-1}/4}\int_{\R^2\setminus B(0,a_v^{-1}/2)}e^{-|y|}\,dy.
\end{align*}
In the last inequality, we have used that $|x-y|\geq a_v^{-1}/4$ in such a way that $|x-y|\geq 2$ for $v$ large enough so as to make $(1+2\ln v)^\alpha/4\geq 2$ and the use $g_m(y)\leq Ce^{-|y|}$ for $|y|\geq 2$.
In both cases ($|x|\leq a_v^{-1}/4$ or $|x|> a_v^{-1}/4$) and since $\alpha>1$, we can choose $\beta>1$ such that  $\sup_{x\in\R^2}|k_m(x)-k(v,x/v)|\leq C  v^{-\beta}$ for some (possibly different) constant $C$. Observe now that similar bounds hold for the gradients $\partial k_m$ and $\partial g_m$ so that we can reproduce the same argument and get for some constant $C$
\begin{equation}
\sup _{x\in\R^2}|\partial k_m(vx)-\partial k(v,\cdot/v)|\leq C v^{-\beta}.
\end{equation} 
[B.2] easily follows from the above uniform convergence; the power decay of the supremum here is not necessary, but we shall use it later. The Seneta-Heyde renormalization thus holds for the family $(X^{sc}_t)_t$. Let us denote by $
M^{\gamma,sc}_t$ and $M^{',sc}_t$ the random measures associated with the field $(X_t^{sc})_t$, sticking to these notations for the rest of the paper.

Let us treat now the long-range correlations $(X^{lc}_t)_t$. The covariance kernel matches:
$$K^{lc}_t(x)=\E[X^{lc}_t(0)X^{lc}_t(x)]=\int_1^{e^t}\int_{\R^2}g_m(y-xv)g_m(y)(1-\varphi) (a_v |y-xv| )(1-\varphi) (a_s |y| )dy\,\frac{dv}{v}$$
and we can set
$$h(v,x)= \int_{\R^2}g_m(y-xv)g_m(y)(1-\varphi) (a_v |y-xv| )(1-\varphi) (a_s |y| )dy,$$ in such a way that
$K^{lc}_t(x)=\int_1^{e^t} \frac{h(v,x)}{v}\,dv.$ We can proceed as we did above to show that
$$\sup _{x\in\R^2}|\partial h(v,x/v)|+|h(v,x/v)|\leq  C v^{-\beta}.$$
Now we use the fact that, for each fixed $x$, the process $t\mapsto (X^{lc}_t(x)-X^{lc}_t(0))_t$ is a martingale to deduce via the Doob inequality that for $q\geq 2$:
\begin{align*}
\E[\sup_t(X_t^{lc}(x)-X^{lc}(0))^q]\leq &C_q\Big(\sup_{t\geq 0}\int_1^{e^t}\frac{h(v,0)-h(v,x)}{v}\,dv\Big)^{q/2}\\
\leq &C|x|^{q/2}\Big(\int_1^{\infty}v^{-\beta}\,dv\Big)^{q/2}.
\end{align*}
The Kolmogorov criterion ensures that, almost surely, the family $(X_t^{lc})_t$ is relatively compact in the space of continuous function equipped with the topology of the uniform convergence over compact subsets of $\R^2$. Furthermore, for each rational points the process $(X^{lc}_t(x))_t$ is a martingale uniformly bounded in $L^2$ and therefore converges almost surely. Both of these statements imply that almost surely, the family $(X_t^{lc})_t$ converges uniformly over compact subsets towards a limiting process, which is nothing but 
$$X^{lc}_\infty(x)= \int_1^{\infty} \int_{\R^2} \frac{g_m(y-xv)}{\sqrt{v}}(1-\varphi) (a_v |y-xv| )W(dv,dy).$$

Now we just have to gather the different pieces of the puzzle to conclude. We have proved that for each measurable set $A$, the family $(\sqrt{t}M^{2,sc}_t(A))$ converges in probability towards   $M^{',sc}(A)$. Therefore, for each subsequence, we can extract a subsequence, still indexed with $t$, such that the family of random measures $(\sqrt{t}M^{2,sc}_t(dx))$ almost surely weakly converges towards $M^{',sc}(dx)$.
Since the family  $(X_t^{lc})_t$ converges uniformly over compact subsets towards $X^{lc}_\infty$ we deduce that for each compactly supported continuous $\phi$ we have the almost sure convergence as $t\to \infty$
$$\sqrt{t}\int_{\R^2}e^{2X^{lc}_t(x)-2\E[X^{lc}_t(x)^2]}\phi(x) M^{\gamma,sc}_t(dx)\to \int_{\R^2}e^{2X^{lc}_\infty(x)-2\E[X^{lc}_\infty(x)^2]}\phi(x) M^{',sc}(dx).$$
Observe that this latter quantity is nothing but $\int_{\R^2} \phi(x) M^{'}(dx)$, where $ M^{'}$ is the derivative multiplicative chaos associated with the family $(X_t)_t$. This comes from the fact that the measure $M^{2,sc}(dx)$ identically vanishes. Hence the result.\qed

\begin{remark}
Observe that we have only use the following properties of $k_m$:\\
- there exist some constants $C>0$ and $\alpha >1$ such that $k_m(x)+|\partial k_m(x)|\leq Ce^{-  |x|^{1/\alpha}}$,\\
-$k_m$ admits a nonnegative convolution square root $g_m$ such that $g,\partial g$ are integrable, and $$\sup_{|x|\geq 2}g(x)+|\partial g(x)|<+\infty\quad\quad \text{ and }\quad \quad \int_1^{\infty}v^{-1}\int_{B(0,\ln^\alpha v)}g^2(y)\,dy\,dv<\infty.$$
Therefore Theorem \ref{seneta} holds for star-scale invariant kernels satisfying these assumptions.
\end{remark}
  
\vspace{2mm} 
 {\bf Application to the cut-off family with covariance \eqref{covMFF2}.} At first sight, the situation may appear different from the previous case because the structure of $p(t,x,y)$ is seemingly different from $k_m(u(x-y))/u$. Let us make the change of variables $u=v^{-2}$ in \eqref{covMFF2}   to get:
\begin{align*} 
G_{t}(x-y)=\int_{0}^{e^{t}} e^{-\frac{|x-y|^2v^2}{2}}\frac{e^{-\frac{m^2}{2v^2}}}{v}\,dv
\end{align*}
This expression  takes on the form $G_{t}(x)=\int_{0}^{e^{t}}\frac{k(vx)}{v}e^{-\frac{m^2}{2v^2}}\,dv$ with $k(x):=e^{-\frac{|x|^2}{2}} $. We can get rid of the smooth part $\int_{0}^{1}\frac{k(vx)}{v}e^{-\frac{m^2}{2v^2}}\,dv$ (for instance by putting it in the long range correlation part) so that we only have to deal with $\int_{1}^{e^{t}}\frac{k(vx)}{v}e^{-\frac{m^2}{2v^2}}\,dv$.  Observe also that the square root convolution of $k$ is still a Gaussian density. It is now clear that we can repeat the above arguments up to the fact that we deal with $e^{-\frac{m^2}{2v^2}}$ instead of $1$ but this is harmless as only the behaviour for large $v$  matters.\qed

\vspace{3mm}

\noindent {\it Proof of Theorem \ref{senetapower}.}
This context is quite close to that we considered in the proof of Theorem \ref{seneta}. The proof remains unchanged, provided that we make the following modifications. First remark that the variance of $X_t$ is no longer $t$, but $K_t(0)$: this is not a real problem, as $K_t(0)$ satisfies $t-c\leq K_t(0)\leq t+c$ for some constant $c>0$, thanks to Assumption [B.1]. In particular, $Y^\beta$ is no longer a Bessel process but rather a time changed Bessel process with a change of time close to $t$ (within two constants). Then follow the proof without change until Lemma  \ref{lem:main}. Then replace the definitions of $R^{\beta,i}_\varepsilon(y),R^{\beta,c}_\varepsilon(y)$ and $Z^{\beta,i}_t(x),Z^{\beta,c}_t(x)$ in the proof of this lemma  by
\begin{align*}
R^{\beta,i}_t(x)&=\int_{B(x,\psi(e^{h_t}))}\one_{\{\tau^\beta_z>t\}}e^{2X_t(z)-2K_t(0)}\,dz & R^{\beta,c}_t(x)&=R^\beta_t(A)-R^{\beta,i}_t(x)  \\
Z^{\beta,i}_t(x)&=\int_{B(x,\psi(e^{h_t})) }f_t^\beta(z)\,dz&Z^{\beta,c}_t(x)&=Z^\beta_t(A)-Z^{\beta,i}_t(x).
\end{align*}
The definitions of the sets $E^1_t, E^2_t, E^3_t, E^4_t, E^5_t$ and $E_t$ remain unchanged. A few lines below, observe that the sigma-algebras $\mathcal{G}^{out,a}_{h_t}$ and $\mathcal{G}^{in}_{h_t}$ are now independent
for $a>\psi(e^{h_t})$. But, conditionally on $\mathcal{G}_{h_t}$, $\frac{R^{\beta,c}_t(x)}{Z^{\beta,c}_t(x)}$ still remains $\mathcal{G}^{out,\psi(e^{h_t})}_{h_t}$-measurable so that we can then follow the proof until Lemma  \ref{deadlycomput2}. In this lemma, replace the definition of $E^6_t$ by 
$$E^6_t=\{h_t^{1/7}\leq \inf_{u\in[h_t,a_t]}Y_u^\beta(x)\} \quad \text{ with }\quad a_t=\ln \psi^{-1}(e^{-t}).$$
With this new definition of $E^6_t$, the only part that needs to be adapted is the treatment of $A^1_t$ in the proof of Lemma \ref{deadlycomput2} but this is painless. To do this, split once again $Z^{\beta,i}_t(x)$ in two parts:
$$ Z^{\beta,i,1}_t(x)=\int_{B(x,e^{-t})\cap A}f_t^\beta(z)\,dz,\quad \text{ and }\quad Z^{\beta,i,2}_t(x)=\int_{C(x,e^{-t},\psi(e^{h_t}))\cap A}f_t^\beta(z)\,dz,$$ where $C(x,e^{-t},\psi(e^{h_t}))$ still stands for the annulus $B(x,\psi(e^{h_t}))\setminus B(x,e^{-t})$. Then the treatment of $Z^{\beta,i,1}_\varepsilon$ is the same and for $Z^{\beta,i,2}_t$, replace everywhere $|x-w|^{-1}$ by $\psi^{-1}(|x-w|)$. At the very end of the proof of this lemma, you get (for some irrelevant constant $C'$, which may change along lines):
 \begin{align*}
\E_{\Theta^\beta_t}&[Z_t^{\beta,i,2}(x)\one_{E^6_t}] \nonumber\\  
\leq& \frac{C'e^{-\sqrt{2d}h_t^{1/7}}}{  |A|}\int_{ A} \int_{C(x,e^{-t},\psi(e^{h_t}))} \big(\psi^{-1}(|x-w|)\big)^d(1+\ln\psi^{-1}(|x-w|))\,dwdx \\
\leq & C'e^{-\sqrt{2d}h_t^{1/7}} \int_{e^{-t}}^{ \psi(e^{h_t})}(\psi^{-1}(r))^d (1+\ln\psi^{-1}(r)  )  r^{d-1} \,dr \\
\leq & C'e^{-\sqrt{2d}h_t^{1/7}} \int^{\psi^{-1}(e^{-t})}_{ e^{h_t}}u^d (1+\ln u )  \psi(u)^{d-1} \psi'(u)\,du.
\end{align*}
Now use the fact that $-\psi'(u)\leq C'(1+\ln u)^{\alpha}u^{-2}$ to get:
 \begin{align*}
\E_{\Theta^\beta_t} [Z_t^{\beta,i,2}(x)\one_{E^6_t}]\leq &       C'e^{-\sqrt{2d}h_t^{1/7}}  \int^{\psi^{-1}(e^{-t})}_{ e^{h_t}}\frac{(1+\ln u)^{\alpha(d+1)}}{u}\,du 
\leq   C'e^{-\sqrt{2d}h_t^{1/7}} t^{\alpha(d+1)+1}.
\end{align*}
Since this quantity decays towards $0$ faster than any power of $t$, we can conclude the whole proof similarly.\qed

\subsection{Gaussian Free Field}
Recall that the normalized Green function over the domain $D$ with Dirichlet boundary conditions is given by  \eqref{bounGreen}. We consider   a white noise $W$ distributed on $\R^2\times \R_+$ and define:
 $$X_t(x)=\sqrt{\pi}\int_{D}\int_{e^{-2t}}^{+\infty}p_D(\frac{s}{2},x,y)W(dy,ds).$$
 
The first observation we make is that the fluctuations of the GFF along the boundary of $D$ are rather tricky to control. Yet this is not a problem: to construct the measure $M'$ on $D$ we just need to construct it on every domain $D'\subset D$ such that ${\rm dist}(D',D^c)>0$ and then consider an increasing sequence $(D_n)_n$ of such domains such that $\bigcup_nD_n=D$. This is a rather standard argument, so that we will not discuss it further.
 
In what follows, we fix a subset $D'\subset D$ such that ${\rm dist}(D',D^c)=\delta>0$. Concerning the construction of the derivative multiplicative chaos, the proof is the same as that explained for the MFF: just use the 4 Brownian motion decomposition to handle the long range correlations. One may also ask how to prove the second statement of \cite[Proposition 19]{DRSV} as the variance of the family $(X_t(x))_t$ depends on the point $x$. This is not that complicated: for each $x\in D'$, define $F(t,x)=\pi\int_{e^{-2t}}^{\infty} p_D(s,x,x)\,ds$. Observe that, on $D'$, we can find a constant $C$ such that for all $x\in D'$ and $t\geq 0$, $t-C\leq F(t,x)\leq t+C$. Since the mapping $t\mapsto F(t,x)$ is continuous and strictly increasing, we can consider the inverse mapping $F^{-1}(t,x)$, which is also close to $t$ within two constants. Then the Gaussian process $Y_t(x)=X_{F^{-1}(t,x)}(x)$ has constant variance $t$. We can then use the comparison argument in \cite{DRSV} with a cone constructed log-correlated Gaussian field (based on \cite[Lemma 18]{DRSV}) to prove Proposition 19. 
 
Concerning the Seneta-Heyde renormalization, the strategy is quite similar to the MFF. We first get rid of the integral between $1$ and $\infty$  and consider
$$X^1_t(x)=\sqrt{\pi}\int_{D\times [e^{-2t},1[}p_D(\frac{s}{2},x,y)W(dy,ds).$$
Then we define another Gaussian field
$$X^{{\rm MFF}}_t(x)=\sqrt{\pi}\int_{\R^2\times [e^{-2t},1[}p(\frac{s}{2},x,y)W(dy,ds).$$
Observe that this is the process we considered in the case of the MFF since we consider here the transition density of the whole-plane planar Brownian motion (up to the irrelevant term $e^{-\frac{m^2}{4}s}$ term). 
In order to conclude as above that the Seneta-Heyde norming holds for the GFF, we have to show that the difference $(X^1_t-X^{{\rm MFF}}_t)_t$  converges uniformly on $D'$ towards a limiting continuous Gaussian process. In fact, it is sufficient to show that $(\bar{X}^1_t-\bar{X}^{{\rm MFF}}_t)_t$ converges uniformly on $D'$ where we set 
$$\bar{X}^1_t(x)=\sqrt{\pi}\int_{D'\times [e^{-2t},1[}p_D(\frac{s}{2},x,y)W(dy,ds).$$
and 
$$\bar{X}^{{\rm MFF}}_t(x)=\sqrt{\pi}\int_{D'\times [e^{-2t},1[}p(\frac{s}{2},x,y)W(dy,ds).$$
Recall the standard formula \cite[section 3.3]{Morters}
\begin{equation*}
\Delta(s,x,y):=p(s,x,y)-p_D(s,x,y)= \E^{x} \big[ 1_{\{T^x_D  \leq s\}}  p(s-T^{x}_D, B^{x}_{T_D^{x}},y )   \big]
\end{equation*}
where $B^x_{t}$ is a standard Brownian motion starting from $x$ and $T^x_D=\text{inf} \{t \geq 0, \; B^x_t\not \in D \}$. Now, take two points $x$ and $x'$. One can couple two Brownian motions from these two starting points by a standard procedure. The two Brownian motions $B^x_t, B^{x'}_t$ run independently until the first components meet; when these components meet, they follow the same $1d$-Brownian motion. After this, the second components evolve independently until they meet and then both $2d$-Brownian motions $B^x_t$ and $B^{x'}_t$ are the same (see \cite[Lemma 2.15]{GRV} for instance). We denote $\P^{x,x'}$ the corresponding probability measure and $\tau^{x,x'}$ the time where both Brownian motions are the same. Now, we have by definition
\begin{align}    
& | \Delta(s,x',y)-\Delta(s,x,y )| \nonumber\\
& = |\E^{x,x'}  [ 1_{\{T^{x'}_D  \leq s\}} p(s-T^{x'}_D, B^{x'}_{T_D^{x'}},y ) - 1_{\{T^x_D  \leq s\}}  p(s-T^{x}_D, B^{x}_{T_D^{x}},y ) ]  |  \nonumber \\
& = | \E^{x,x'} [ 1_{\{ \tau^{x,x'} >  T^{x}_D \wedge T^{x'}_D\}} 1_{\{T^{x'}_D  \leq s\}} p(s-T^{x'}_D, B^{x'}_{T_D^{x'}},y ) - 1_{ \{\tau^{x,x'} >  T^{x}_D \wedge T^{x'}_D\}} 1_{\{T^x_D  \leq s\}}  p(s-T^{x}_D, B^{x}_{T_D^{x}},y )  ]   |  \nonumber\\
& \leq | \E^{x,x'} [ 1_{\{ \tau^{x,x'} >  T^{x}_D \wedge T^{x'}_D\}} 1_{\{T^{x'}_D  \leq s\}} p(s-T^{x'}_D, B^{x'}_{T_D^{x'}},y ) +1_{ \{\tau^{x,x'} >  T^{x}_D \wedge T^{x'}_D\}} 1_{\{T^x_D  \leq s\}}  p(s-T^{x}_D, B^{x}_{T_D^{x}},y )  ]   |  \nonumber\\
& \leq \frac{e^{-\delta^2/2s}}{\pi s}  \P^{x,x'} ( \tau^{x,x'} >  T^{x}_D \wedge T^{x'}_D  )\label{yiyi}
\end{align}
An easy calculation shows that the law $\frac{\tau^{x,x'}}{|x-x'|^2}$ is independent of $x,x'$ and matches that of a stable law with index stability $1/2$, call it $S$. Then, by using the Markov inequality, we deduce, for $\beta\in ]0,2[$, $ \alpha \in]0,1/2[$ and $\chi>0$, 
\begin{align}
 \P^{x,x'} ( \tau^{x,x'} >  T^{x}_D \wedge T^{x'}_D  )\leq & \P^{x,x'} ( \tau^{x,x'} >  |x-x'|^\beta )+ \P^{x} ( T^{x}_D  \leq   |x-x'|^\beta )+  \P^{x'} ( T^{x'}_D\leq   |x-x'|^\beta ) \nonumber\\
 \leq & |x-x'|^{(2-\beta)\alpha}\E[S^\alpha]+2|x-x'|^{\beta\chi}\sup_{x\in D'}\E[(T^x_D)^{-\chi}].
 \label{yiyi2}
\end{align}
If $\delta>0$ stands for ${\rm dist}(D',C^c)>0$, we have  
\begin{equation}
\sup_{x\in D'}\E[(T^x_D)^{-\chi}]\leq \sup_{x\in D'}\E[(T^x_{B(x,\delta)})^{-\chi}]=\E[(T^0_{B(0,\delta)})^{-\chi}]<+\infty. \label{yiyi3}
\end{equation}
The latter quantity is finite for all $\chi>0$ (standard argument: compute the law of $T^0_{B(0,\delta)}$ by applying the stopping time theorem to the exponential martingale of the Brownian motion).
By combining \eqref{yiyi}+\eqref{yiyi2}+\eqref{yiyi3}, one concludes as explained above by the Kolmogorov criterion that $(\bar{X}^1_t - \bar{X}^{{\rm MFF}}_t)_t$ converges uniformly on $D'$ towards a limiting continuous Gaussian process.   \qed

\vspace{2mm}
\noindent {\it Proof of Theorem \ref{conformal}.}
Given another planar domain $\widetilde{D}$ and a conformal map $\psi:\widetilde{D}\to D$, we will denote by $(X^\psi_t)_{t\geq 0}$ the random field  defined by 
$$X^\psi_t(x)=X_t(\psi(x)).$$
This family is a white noise approximation of the centered GFF $X\circ \psi$  defined on $\widetilde{D}$. Then  for any function $\varphi$ continuous and compactly supported  on $\widetilde{D}$, we consider:
\begin{align*}
&M^{X\circ \psi +2\ln|\psi'(x)|,\widetilde{D}}_t(\varphi)\\
&=\int_{\widetilde{D}}\varphi(x) (2\E[X^\psi_t(x)^2]-X^\psi_t(x)-2\ln |\psi'(x)|) e^{2 X^\psi_t(x) +4\ln|\psi'(x)|-2\E[X^\psi_t(x)^2]}C(x,\tilde{D})^{2} \,dx.
\end{align*}
When $\psi$ is the identity we simply write $M^{X,D}_t$.
This defines a martingale with respect to the parameter $t$. By using the standard rule, for $x\in\tilde{D}$, $| C( \psi(x),D)  |= |\psi'(x)| | C(x,\tilde{D})  |$ where  $|\psi'(x)|^2$  is the Jacobian  of the mapping $\psi:\widetilde{D}\to D$, we get  the relation $$M^{X\circ \psi +2\ln|\psi'|,\widetilde{D}}_t(\varphi)=M^{X,D}_t(\varphi\circ \psi^{-1}) -2 \int_D\ln |\psi'\circ\psi^{-1}(x)|\, \varphi\circ \psi^{-1}(x)e^{2X_t(x)-2\E[X_t(x)^2]}C(x,D)^2\,dx.$$
By using the fact that the random measure $e^{2X_t(x)-2\E[X_t(x)^2]}C(x,D)^2\,dx$  almost surely weakly converges towards $0$ as $t\to\infty$, we can pass  to the limit as $t\to \infty$ in the above relation and get the almost sure weak convergence of the family of random measures $(M^{X\circ \psi+2\ln|\psi'| ,\widetilde{D}}_t(dx))_t$ towards a limiting measure $M^{X\circ \psi+2\ln|\psi'| ,\widetilde{D}}(dx)$ satisfying
$$M^{X\circ \psi+2\ln|\psi'| ,\widetilde{D}}(\varphi)=M^{X  ,D}(\varphi\circ\psi^{-1}) $$ for every continuous and compactly supported function $\varphi$ on $\widetilde{D}$.\qed


\hspace{10 cm}


\begin{thebibliography}{20}
\bibitem{AidShi}
A\"{\i}d\'ekon E., Shi Z.: The Seneta-Heyde scaling for the branching random walk, arXiv:1102.0217v2.

\bibitem{allez}
Allez R., Rhodes R., Vargas V.: Lognormal $\star$-scale invariant random measures, {\it Probability Theory and Related Fields}, to appear, arXiv:1102.1895v1.

\bibitem{Al} Alvarez-Gaum\'e L., Barb\'on, J. L. F., and  Crnkovi\'c \v{C}.: A proposal for strings at $D>1$,  \emph{Nucl.
Phys.} \textbf{B394}, 383 (1993).

\bibitem{Amb} Ambj{\o}rn J., Durhuus B., and Jonsson T.: A solvable 2d gravity model with $\gamma >0$,  \emph{Mod. Phys. Lett.}
A 9, 1221 (1994).

\bibitem{bacry}
Bacry E., Muzy J.F.: Log-infinitely divisible multifractal processes, \emph{Comm. Math. Phys.}, \textbf{236} (2003) no.3, 449-475.

\bibitem{BJRV}
Barral J., Jin X., Rhodes R., Vargas V.: Gaussian multiplicative chaos and KPZ duality, arXiv:1202.5296v2.

\bibitem{Bar}
Barral, J., Mandelbrot, B.B.: Multifractal products of cylindrical pulses, \emph{Probab. Theory
Relat. Fields} \textbf{124} (2002), 409-430.


\bibitem{BK}
Barral J., Kupiainen A., Nikula M., Saksman E., Webb C.: Critical Mandelbrot's cascades, arXiv:1206.5444v1.









\bibitem{Benj}
Benjamini, I., Schramm, O.: KPZ in one dimensional random geometry of multiplicative cascades, Communications in Mathematical Physics, vol. 289, no 2, 653-662, 2009.



\bibitem{Louisdor}
Biskup M., Louidor O.: Extreme local extrema of the two-dimensional discrete Gaussian free field,  	arXiv:1306.2602.


\bibitem{Bram}
Bramson M., Ding J., Zeitouni O.: Convergence in law of the maximum of the two-dimensional discrete Gaussian free field, arXiv:1301.6669.

\bibitem{BKZ}
Br\'ezin E., Kazakov V.A. ,  Zamolodchikov Al.B.: Scaling violation in a field theory of closed strings in one physical dimension, \emph{Nuclear Physics} \textbf{B338},  673-688 (1990).




\bibitem{chelkak}
Chelkak D., Smirnov S.: Discrete complex analysis on isoradial graphs, Advances in Mathematics 228 (2011) 1590Ð1630. 



\bibitem{Das}  Das S. R., Dhar A., Sengupta A. M., and Wadia S. R.: New critical behavior in $d=0$ large-$N$ matrix models, \emph{Mod.
Phys. Lett.  A} \textbf{5}, 1041 (1990).

\bibitem{cf:Da} David, F.: Conformal Field Theories Coupled to 2-D Gravity in the Conformal Gauge, \emph{Mod. Phys. Lett. A}, \textbf{3} (1988).

\bibitem{DFGZ} Di Francesco P., Ginsparg P., Zinn-Justin J.: 2D gravity and random matrices, \emph{Physics Reports} \textbf{254}, p. 1-133 (1995).

\bibitem{DistKa} Distler J.,   Kawai H.: Conformal Field Theory and 2-D Quantum Gravity or Who's Afraid of Joseph Liouville?, \emph{Nucl. Phys.} \textbf{B321} 509-517 (1989).

\bibitem{DRSV}
Duplantier B., Rhodes R., Sheffield S., Vargas V.: Critical Gaussian Multiplicative Chaos: convergence of the derivative martingale, arXiv:1206.1671.



\bibitem{Dup:houches}
Duplantier B.: A rigorous perspective on Liouville quantum gravity and KPZ, in \emph{Exact Methods in Low-dimensional Statistical Physics and Quantum Computing},
J. Jacobsen, S. Ouvry, V. Pasquier, D. Serban, and L.F. Cugliandolo, eds., Lecture Notes of the Les Houches Summer School: Volume 89, July 2008, Oxford University Press (Clarendon, Oxford) (2010).

\bibitem{BDMan} Duplantier B.: Conformal fractal geometry and boundary quantum gravity, in \emph{Fractal geometry and applications: a
jubilee of Beno\^it Mandelbrot}, Part 2 (Amer. Math. Soc.,
Providence, RI, 2004), vol. 72 of \emph{Proc. Sympos. Pure
Math.}, pp. 365--482, arXiv:math-ph/0303034.


\bibitem{cf:DuSh} Duplantier, B., Sheffield, S.: Liouville Quantum Gravity and KPZ, \emph{Inventiones Mathematicae}, 2011, 185 (2), 333-393.

\bibitem{PRL} Duplantier, B., Sheffield, S.: Duality and KPZ in Liouville Quantum Gravity, \emph{Physical Review Letters}, {\bf 102}, 150603 (2009).


\bibitem{Dur} Durhuus B.: Multi-spin systems on a randomly triangulated surface, \emph{Nucl. Phys.} \textbf{B426}, 203 (1994).



\bibitem{falc}
 Falconer K.J.: The geometry of fractal sets, Cambridge University Press, 1985.

\bibitem{Fan}
Fan A.H., {\it Sur le chaos de L\'evy d'indice $0<\alpha<1$}, \emph{Ann. Sciences Math. Qu\'ebec}, vol 21 no. 1, 1997, p. 53-66.

\bibitem{GRV}
Garban C., Rhodes R., Vargas V.:  Liouville Brownian Motion, arXiv:1301.2876v2 [math.PR].

\bibitem{GM} Ginsparg P. and Moore G.: Lectures on 2D gravity and 2D string theory, in \textit{Recent direction in particle
theory}, Proceedings of the 1992 TASI, edited by J. Harvey
and J. Polchinski (World Scientific, Singapore, 1993).

\bibitem{GZ} Ginsparg P., Zinn-Justin J.: 2D gravity $+$ 1D matter, \emph{Physics Letters} \textbf{B240},  333-340 (1990).

 \bibitem{GrossKleban}
Gross, D.J. and Klebanov I.R.: One-dimensional string theory on a circle, \emph{Nuclear Physics} \textbf{B344} (1990) 475--498.

\bibitem{GrossM} Gross D.J., Miljkovi\'c N.: A nonperturbative solution of $D = 1$ string theory,
\emph{Physics Letters} \textbf{B238},  217-223 (1990).

\bibitem{GubserKleban}
Gubser S.S., Klebanov I.R.: A modified $c = 1$ matrix model with new critical behavior,
\emph{Physics Letters}  \textbf{B340} (1994) 35--42.

\bibitem{heyde}
Heyde  C.C.: Extension of a result of Seneta for the super-critical Galton-Watson process, \emph{Ann. Math. Statist.} 41, 739-742 (1970).

\bibitem{HuShi}
Hu Y., Shi Z.: Minimal position and critical martingale convergence in branching random walks, and directed polymers on disordered trees, \emph{Annals of Probability}, \textbf{37} (2) (2009) 742-789.

\bibitem{Jain}Jain S. and  Mathur S. D.: World-sheet geometry and baby universes in 2-D quantum gravity, \emph{Phys. Lett.} \textbf{B 286}, 239 (1992).



\bibitem{cf:Kah} Kahane J.-P.: Sur le chaos multiplicatif,
 \emph{Ann. Sci. Math. Qu{\'e}bec}, \textbf{9} no.2 (1985), 105-150.


 \bibitem{KKK} Kazakov V.,  Kostov I., and  Kutasov D.:  A Matrix Model for the 2d Black Hole,
in \textit{Nonperturbative Quantum Effects 2000}, JHEP Proceedings; A matrix model for the two-dimensional black hole, \emph{Nucl. Phys.} \textbf{B622}  (2002) 141-188.

\bibitem{Kle}
Klebanov, I.: String theory in two dimensions. arXiv:hep-th/9108019.

\bibitem{Kleb1}
Klebanov I.R.: Touching random surfaces and Liouville gravity, \emph{Phys. Rev. D 51}, 1836Ð1841 (1995).

\bibitem{Kleb2}
Klebanov I.R., Hashimoto A.: Non-perturbative Solution of Matrix Models Modified by Trace-squared Terms, \emph{Nucl. Phys}. \textbf{B434} (1995) 264-282.

\bibitem{Kleb3}
Klebanov I.R., Hashimoto A.: Wormholes, Matrix Models, and Liouville Gravity, \emph{Nucl. Phys.} (Proc. Suppl.) \textbf{45}B,C (1996) 135-148.

\bibitem{cf:KPZ} Knizhnik V.G., Polyakov A.M., Zamolodchikov A.B.: Fractal structure of 2D-quantum gravity, \emph{Modern Phys. Lett A}, \textbf{3}(8) (1988), 819-826.

\bibitem{kostov91} Kostov I.K.: Loop amplitudes for nonrational string theories, \emph{Phys. Lett.}  \textbf{B266}, 317-324 (1991).

\bibitem{kostov92} Kostov I.K.:  Strings with discrete target space, \emph{Nucl. Phys.} \textbf{B376}, 539-598 (1992).

\bibitem{Kostov:houches}
Kostov I.K.: Boundary Loop Models and and 2D
Quantum Gravity, in \emph{Exact Methods in Low-dimensional Statistical Physics and Quantum Computing},
J. Jacobsen, S. Ouvry, V. Pasquier, D. Serban, and L.F. Cugliandolo, eds., Lecture Notes of the Les Houches Summer School: Volume 89, July 2008, Oxford University Press (Clarendon, Oxford) (2010).

\bibitem{KS} Kostov I.K., Staudacher M.: Multicritical phases of the $O(n)$ model on a random lattice, \emph{Nucl. Phys.} \textbf{B384}, 459-483 (1992).

\bibitem{lacoin}
Lacoin H, Rhodes R., Vargas V.: Complex Gaussian multiplicative chaos, arxiv.



\bibitem{madaule}
Madaule T.: Maximum of a log-correlated Gaussian field, arXiv:1307.1365v2 [math.PR].

\bibitem{Mol}
 Molchan, G.M.: Scaling exponents and multifractal dimensions for independent random cascades, \emph{Communications in Mathematical Physics}, \textbf{179} (1996), 681-702.

\bibitem{Morters}
Morters P. Peres Y.: Brownian motion, \emph{Cambridge University press}, (2010). 

\bibitem{Nak}
Nakayama Y.: Liouville Field Theory -- A decade after the revolution, \emph{Int. J. Mod. Phys.} \textbf{A19}, 2771 (2004).

\bibitem{Parisi}
Parisi G.: On the one dimensional discretized string, \emph{Physics Letters} \textbf{B238},  209-212
(1990).

\bibitem{Polch} Polchinski J., Critical behavior of random surfaces in one dimension, \emph{Nuclear Physics} {\bf B346} (1990) 253--263.

\bibitem{sohier} Rhodes R., Sohier J., and Vargas V.: $\star$-scale invariant random measures, to appear in Annals of Probability, arXiv:1201.5219v1.

\bibitem{rhovar} Rhodes R., Vargas, V.: Multidimensional multifractal random measures, \emph{Electronic Journal of Probability}, \textbf{15} (2010), 241-258.

\bibitem{review} 
Rhodes R., Vargas, V.: Gaussian multiplicative chaos and applications: a review, arxiv.

\bibitem{cf:RhoVar} Rhodes, R. Vargas, V.: KPZ formula for log-infinitely divisible multifractal random measures,  \emph{ESAIM Probability and Statistics}, 15, 2011, 358-371.

%

\bibitem{cf:RoVa} Robert, R., Vargas, V.: Gaussian Multiplicative Chaos revisited, \emph{Annals of Probability}, \textbf{38} 2 (2010), 605-631.

\bibitem{SchShe}
Schramm O. Sheffield S.: A contour line of the continuum Gaussian free field, \emph{Probability Theory and related Fields} (2012).

\bibitem{seneta}
Seneta E.: On recent theorems concerning the supercritical Galton-Watson process, \emph{Ann. Math. Statist.}, 39, 2098-2102, 1968.

\bibitem{She07}
Sheffield S.: Gaussian free fields for mathematicians, \emph{Probab. Th. Rel. Fields}, 139:521--541, 2007.


\bibitem{sugino}
Sugino F., Tsuchiya O.:  Critical behavior in $c=1$ matrix model with branching interactions, \emph{Mod. Phys. Lett.} \textbf{A9} (1994) 3149-3162.
\end{thebibliography}
\end{document}